%% file: agt-2-7.tex
\documentclass{gtart}

\input agtout

\lognumber{7}
\volumenumber{2}
\volumeyear{2002}
\papernumber{7}
\published{28 February 2002}
\pagenumbers{137}{155}
\received{1 October 2001}
\accepted{6 February 2002}

\usepackage{amssymb}

\let\Bbb\mathbb

\def\C{{  \Bbb C\,}}

\def\span{{  \rm {span}}}
\def\id{{  \rm {id}}}

  \def\skipaline{\ppar}

  \def\Int{{ \rm  {Int}}}
   \def\Hom{{ \rm  {Hom}}}

  \def\rel {{ \rm  {rel}}}

   \def\exp {{ \rm  {exp}}}

   \def\dim  {{ \rm  {dim}}}

\def\Coker {{ \rm  {Coker}}}
   \def\Tors {{ \rm  {Tors}}}
   \def\rk {{ \rm  {rk}}}
\def\id {{ \rm  {id}}}

   \def\proj{{ \rm  {pr}}}

  \def\Ker {{ \rm  {Ker}}}

  \def\det {{ \rm  {det}}}

   \def\aug {{ \rm  {aug}}}

\begin{document} 

\title{A norm for the cohomology of 2-complexes}

\author{Vladimir Turaev}
\address{IRMA, Universit\'e Louis Pasteur  -- CNRS\\7 rue Ren\'e Descartes,
67084 Strasbourg, France}

\asciiaddress{IRMA, Universite Louis Pasteur -- CNRS\\7 rue Rene Descartes,
67084 Strasbourg, France}
\email{turaev@math.u-strasbg.fr}

\begin{abstract}   We  introduce a   norm on the real 1-cohomology of finite
2-complexes determined by the Euler characteristics of graphs on these
complexes. We also introduce
twisted Alexander-Fox polynomials of groups and  show that they give
rise to
norms on the real 1-cohomology of   groups. Our main theorem states
that for a
finite 2-complex $X$,
the
   norm on   $H^1(X;\Bbb R)$   determined
by
graphs on $X$ majorates the
Alexander-Fox norms    derived from $\pi_1(X)$.
  \end{abstract}

\asciiabstract{We introduce a norm on the real 1-cohomology of finite
2-complexes determined by the Euler characteristics of graphs on these
complexes. We also introduce twisted Alexander-Fox polynomials of
groups and show that they give rise to norms on the real 1-cohomology
of groups. Our main theorem states that for a finite 2-complex X, the
norm on H^1(X; R) determined by graphs on X majorates the
Alexander-Fox norms derived from \pi_1(X).}

\maketitle

\section*{Introduction}

  We  introduce a (possibly degenerate) norm on the real 1-cohomology of
finite
2-complexes. The definition of this norm is  similar to Thurston's
definition of
a norm on the 2-homology of 3-manifolds. The key difference is that
instead of
surfaces in 3-manifolds we consider graphs on 2-complexes.
In many instances the resulting theory is similar to but simpler than
the one of
Thurston.

In generalization of the standard Alexander-Fox polynomial of groups we
introduce
twisted Alexander-Fox polynomials. We show that they determine
norms on the real 1-cohomology of   groups.

Our main result   is a comparison theorem which states that for a
finite
2-complex $X$,
the
   norm on   $H^1(X;\Bbb R)$   determined
by
graphs on $X$ majorates the
Alexander-Fox norms    derived from $\pi_1(X)$.

This result is a cousin of the classical Seifert inequality in knot
theory  which says that the genus of a knot $K\subset S^3$ is greater
than or
equal
to the half of the span of the Alexander  polynomial of $K$. A more
general
estimate
from
below for the Thurston norm
appeared in the Seiberg-Witten theory in dimension 3,
see
\cite {Au},  \cite {Kr}, \cite {KM1},
\cite {KM2}.  This estimate is
a
3-dimensional version of the much deeper adjunction inequality in
dimension 4.  A related (weaker) result  in dimension 3 appeared also
in \cite {McM}.

We state here a sample application of our main
theorem to codimension 1 submanifolds of triangulated manifolds.
Let  $M$ be a  closed
connected oriented
triangulated manifold    of dimension $ m\geq 3$.   Let  $S\subset M$
be a
closed oriented
  $(m-1)$-dimensional submanifold of $M$   representing a non-zero
element $s\in
H_{m-1}(M;\Bbb Z)=
H^1(M;\Bbb Z)$. Let $n$ be the maximal positive integer dividing $s$ in
$
H^1(M;\Bbb Z)$. Assume that $S$ intersects the 2-skeleton $M^{(2)}$
transversely
along a (finite 1-dimensional) CW-space $\Gamma=S\cap M^{(2)}$.  If
$\pi_1(M)=\pi_1(S^3\backslash K)$ where $K$ is a knot in $S^3$ then
   $\vert \chi (\Gamma) \vert\geq n(d-1)$ where $d$ is the span of the
Alexander
polynomial of $K$.  For example, if $\pi_1(M)=\langle x, y: x^p
y^q=1\rangle$  is the group of a torus $(p,q)$-knot with relatively
prime
integers
$p,q\geq 2$ then $\vert \chi (\Gamma) \vert\geq n(pq-p-q)$.

\section{A norm on the 1-cohomology of a 2-complex}

{\bf 1.1\qua Two-complexes}\qua
  By a {\sl graph} we    mean a finite   CW-complex of dimension $\leq
1$.
By a {\sl finite 2-complex} we  mean the underlying topological space
of a
finite
2-dimensional
CW-complex   such
that
each its point     has a neighborhood
homeomorphic to  the cone over a graph. The   latter condition is aimed
at
eliminating
all kinds of local wilderness.
Examples of finite
2-complexes: compact surfaces;    2-skeletons of   finite simplicial
spaces;
products of   graphs  with a closed interval.

We define two subspaces $\Int X$ and $\partial X$ of a finite 2-complex
$X$.
The subspace $\Int X\subset X$   consists of the
points
which have a neighborhood homeomorphic to $\Bbb R^2$. Clearly,  $\Int
X$ is a
2-manifold
with finite number of components. Its complement   $X\backslash \Int X$
is a
graph contained in the 1-skeleton of any CW-decomposition of $X$.

The {\sl boundary} $\partial X$ of $X$ is the closure in $X$ of the set
of
all
points of $X\backslash \Int X$ which have  an open neighborhood   in
$X$
homeomorphic to $\Bbb R$ or
to   $\Bbb
R^2_+=\{(a,b)\in \Bbb R^2, b\geq 0\}$.    A simple local analysis shows
that
$\partial X$ is a graph contained in the 1-skeleton of any
CW-decomposition of
$X$.
  If
$X$ is
a compact surface
then
$\partial X$ is its boundary in the usual sense and $\Int X=
X\backslash
\partial
X$.

\medskip

{\bf 1.2\qua Graphs on 2-complexes}\qua
     A    graph $\Gamma$ embedded in a finite 2-complex  $X$ is   {\sl
regular}
if   $\Gamma\subset X\backslash \partial X $ and  there are a closed
neighborhood $U$ of
$\Gamma$
in $X\backslash \partial X$ and a
homeomorphism
$U\approx
\Gamma\times [-1,1]$
sending any point $x\in \Gamma $ to $x\times 0$. If $\Gamma$ is
connected then
  $U\backslash \Gamma$ has two components. A choice of one of them
is
called a {\sl coorientation} of $\Gamma$. If $\Gamma$ is not connected
then a {\sl coorientation}
   of $\Gamma$ is a choice of coorientation for all
components
of
$\Gamma$.

   Any vertex of a regular graph
$\Gamma\subset X$ is incident
to
at least two edges of $\Gamma$ (counting with multiplicity). Hence
$\chi
(\Gamma)\leq 0$. Set
  $ \chi_- (\Gamma)=-\chi(\Gamma)\geq 0$.

A cooriented regular graph $\Gamma\subset X$ determines a
1-dimensional cohomology
class
$s_\Gamma\in H^1(X,\partial X)=H^1(X,\partial X;\Bbb Z)$
  as follows. Choose a neighborhood $U$ of
$\Gamma$
and
a
homeomorphism $f:U\to \Gamma\times [-1,1]$ as above so that the
coorientation
of
$\Gamma$ is determined by the components of $U\backslash \Gamma$ lying
   in  $f^{-1}(\Gamma \times (0,1])$.   We define a
map
$g:X/\partial X\to S^1=\{z\in \C,\, \vert z\vert =1\}$ by $g
(X\backslash
U)=-1\in
S^1$
and $g (x)=\exp (\pi i \tilde f(x))$ for $x\in U$ where $\tilde f(x)\in
[-1,1]$ is
the
projection of $f(x)\in \Gamma\times [-1,1]$ to $[-1,1]$. Set  $s_\Gamma
=g^* (s_0)$ where $s_0$ is the   generator of $H^1(S^1)=\Bbb
Z$ determined by the counterclockwise orientation of $S^1$. It is clear
that
$s_\Gamma$ does not depend on the choice of $U$ and $f$. To evaluate
$s_\Gamma$ on the homology class of a path in $X$ whose endpoints
either
coincide or
lie in $\partial X$, one should count
the
algebraic
number of intersections of this path with $\Gamma$. If
$\Gamma=\emptyset$, then
$s_\Gamma=0$.

A simple transversality argument shows that for any     $s\in
H^1(X,\partial X)$ there is  a cooriented regular graph
$\Gamma\subset X$ such that $s=s_\Gamma$. It can be
constructed as
follows. First, one
realises
$s$ as $g^* (s_0)$ for a certain map
$g:X\to S^1$ sending $\partial X$ to $-1\in S^1$. Secondly,  one fixes
a
CW-decomposition of $X$ and deforms $g$ so
that it
maps the 0-skeleton
$X^{(0)}$ of
$X$
into
$S^1\backslash \{1\}$. Then one deforms $g (\rel \, X^{(0)})$ so that
its
restriction
to
the 1-skeleton $X^{(1)}$ of $X$ becomes transversal to the point $1\in
S^1$.
Finally, one deforms $g (\rel \, X^{(1)})$ so that its restriction to
any
2-cell of
$X$
becomes transversal to   $1\in S^1$. Then $\Gamma=g^{-1} (1)$  is a
regular
graph
on $X$ and   $g$ determines its coorientation such that
$s_\Gamma=s$.

\ppar

{\bf 1.3\qua A norm on $H^1(X, \partial X;\Bbb
R)$}\qua
   By a {\sl norm} on a
real
vector
space $V$ we   mean an $\Bbb R$-valued function  $\vert \vert ...
\vert
\vert$
on
$V$ such that $\vert \vert s  \vert \vert\geq 0$ and
$\vert \vert s+s'  \vert \vert\leq \vert \vert s  \vert \vert+\vert
\vert s'
\vert
\vert$ for any $s,s'\in V$. A norm is allowed to be degenerate, i.e.,
to vanish
on
nonzero vectors. A norm $\vert \vert ...  \vert
\vert$
on
$V$ is   {\sl homogeneous}, if  $\vert \vert ks\vert
\vert =\vert k\vert\,
\vert \vert s\vert \vert $ for any $ k\in \Bbb  R, s \in V$. One
similarly
defines
norms on  lattices, the only difference is that in the
definition of
homogeneity $k\in \Bbb Z$.

Let $X$ be a   finite  2-complex.
For   $s\in H^1(X, \partial X)=H^1(X, \partial X;\Bbb Z)$, set
$$\vert\vert s \vert\vert=\min_{\Gamma, s_\Gamma=s} \,  \chi_- (\Gamma)
$$
where
$\Gamma$
runs
over cooriented regular graphs
   in $ X$ such that $s=s_\Gamma$.
The next lemma shows that  $\vert\vert
...
\vert\vert$ is a
   homogeneous norm on $H^1(X, \partial X)$.
It    extends uniquely to a homogeneous continuous
norm  on  $H^1(X, \partial X;\Bbb  R)$    denoted $\vert\vert ...
\vert\vert_X$ or simply $\vert\vert ...
\vert\vert$.

\skipaline 

{\bf 1.4\qua Lemma}\qua {\sl  $\vert\vert ...
\vert\vert$ is a homogeneous norm on $H^1(X,\partial X)$.}

\proof
We verify that
$\vert\vert s+s' \vert\vert\leq \vert\vert s \vert\vert+ \vert\vert s'
\vert\vert$
for any $s,s'\in H^1(X, \partial X)$. Let $\Gamma, \Gamma'$ be
cooriented
regular
graphs
   in $ X$ such that $s=s_\Gamma, s'=s_{\Gamma'}$. We  slightly deform
$\Gamma$
so
that   $\Gamma\cap \Gamma' \subset \Int X$ and each point
$\gamma\in \Gamma\cap
\Gamma'$ is a transversal   intersection
of
an (open) edge of $\Gamma$ with an (open) edge of $\Gamma'$. A {\sl
smoothing}
of
$\Gamma\cup
\Gamma'$ at     $\gamma $  replaces the  crossing at $\gamma$  by the
$\supset
\subset$-type configuration. There is a
unique
smoothing at   $\gamma$  such that the
coorientations of $\Gamma, \Gamma'$
   induce (locally) a coorientation of the resulting graph.
Applying this smoothing at all points of $\Gamma\cap \Gamma'$ we
transform
$\Gamma\cup \Gamma'$ into a cooriented regular graph, $\Gamma''$,
   in $ X$. It is obvious that $s_{\Gamma''}=s+s'$ and
  $\vert\vert s_{\Gamma''} \vert\vert \leq \chi_- (\Gamma'')= \chi_-
(\Gamma)+\chi_- (\Gamma')$.
   Therefore $\vert\vert s+s' \vert\vert\leq \vert\vert s \vert\vert+
\vert\vert
s'
\vert\vert$.

The homogeneity of  $\vert\vert ...
\vert\vert$ is proven by   the same argument as in \cite {Th}, p.103. The key
point is
that if a cooriented regular graph $\Gamma$ in $X$ represents $ks$ with
integer
$k\geq
1$ and $s\in H^1(X, \partial X )$ then $\Gamma$ splits as a disjoint
union of $k$  graphs representing $s$.
This implies that  $\vert\vert ks \vert\vert\geq k \vert\vert s
\vert\vert$.
The opposite inequality is obvious since for any $\Gamma \subset X$
representing
$s$ a
union of $k$ parallel copies of $\Gamma$
represents $ks$.\endproof

  \skipaline

{\bf 1.5\qua  Properties of
$\vert\vert...\vert\vert_X$}\qua
(1)\qua Replacing everywhere
   embedded   graphs in $X$  by immersed graphs, we obtain the same
norm. (By
an immersed graph we mean a graph  in $X$ which locally  looks like  an
embedded
graph
or like  a   transversal crossing of two embedded arcs in $\Int X$.)
The immersed graphs lead to the same norm  because the smoothing of
an
immersed graph  at all its double points
yields an embedded graph with the same  Euler characteristic.

(2)\qua  It is easy to describe the subset of $ H^1(X, \partial X)$
consisting of
the
vectors  with zero norm. Indeed, for a  regular graph $\Gamma\subset X$
we have
$\chi(\Gamma)=0$ if and only if $\Gamma$ is a closed 1-dimensional
submanifold
of
$\Int X$. Therefore the set of vectors in  $ H^1(X, \partial X)$ with
zero norm
coincides with the set
of   vectors $s_\Gamma$ corresponding to   cooriented closed
1-manifolds
$\Gamma\subset \Int X$.
The  argument   in \cite {Th}, p.105  shows that
the set of vectors in  $ H^1(X, \partial X;\Bbb R)$ with zero norm is
the  $\Bbb
R$-linear
  span
of such $s_\Gamma$.

(3)\qua If all components of $\Int X$ are open 2-discs   or
M\"obius bands
then
the norm $\vert\vert s \vert\vert_X$ vanishes only for $s=0$.
The unit ball $\{s\in  H^1(X, \partial X;\Bbb R), \,
\vert\vert
s\vert\vert \leq 1\}$ is then a compact convex  polytope symmetric in
the origin.
It can be defined by a   system of inequalities
$\vert \beta (s)\vert \leq 1$ where $\beta$ runs over a finite subset
of
$H_1(X, \partial X)$. This follows from general properties of norms
taking
integral values on a lattice of maximal rank, see \cite {Th}, p. 106.

(4)\qua If $p:\tilde X \to X$ is an $n$-sheeted covering with $n \geq 2$
then
  $\partial \tilde X=p^{-1}(\partial X)$ and for any
  $s\in H^1(X, \partial X;\Bbb R)$, we have $\vert \vert p^*(s)\vert
\vert = n \,
\vert
\vert s\vert \vert$
where $p^*$ is the induced homomorphism
$ H^1(X, \partial X;\Bbb R)\to  H^1(\tilde X, \partial \tilde X;\Bbb R)
$.
Indeed, if $\Gamma$ is a  cooriented regular  graph in $X$ representing
$s$ then
the graph $p^{-1}(\Gamma)\subset \tilde X$ represents $p^*(s)$.
Therefore
  $\vert \vert p^*(s)\vert \vert \leq  n\, \vert \vert s\vert \vert$.
On the other hand, if $\Gamma'$ is a cooriented regular   graph in
$\tilde X$
representing
$p^*(s)$ then deforming if necessary $\Gamma'$ we can assume that
$p(\Gamma')$ is
  an immersed graph. Smoothing
it at all crossing points  we obtain a cooriented regular graph
$\Gamma\subset
X$
such that $\chi (\Gamma)=\chi (\Gamma')$ and $s_\Gamma=ns$. Hence
$ n\, \vert \vert s\vert \vert
=\vert \vert ns\vert \vert \leq \vert \vert p^*(s)\vert \vert$.

\skipaline
{\bf 1.6\qua A computation from cocycles}\qua   For any
finite
2-complex $X$ and
   $s\in H^1(X, \partial X ) $, we can compute
  $\vert\vert s \vert\vert_X$   in terms of 1-cocycles on $X$. Fix a
CW-decomposition of $X$ and orient all its edges ($=$ open 1-cells).
Consider  a
$\Bbb
Z$-valued cellular 1-cocycle $k$ on $(X, \partial X)$.
Set $\vert k\vert=\sum_e  (n_e/2 -1) \vert k(e)\vert$
  where $e$ runs over all edges  of $X$ not lying on $\partial X$,
$n_e\geq 2$ is
the
number of 2-cells of $X$ adjacent to $e$ (counted with multiplicity),
and
$k(e)\in \Bbb Z$ is the value of $k$ on $e$.
  We claim that $\vert \vert s\vert\vert_X=\min_k \vert k\vert$ where
$k$ runs
over
all
   cellular 1-cocycles on   $(X, \partial X)$ representing $s$. This
reduces the
computation of $\vert \vert s\vert\vert$ to a   standard
algorithmically
solvable minimization problem
   on a lattice.

  We first prove that $\vert \vert s\vert\vert\leq   \vert k\vert$ for
any $k$ as
above.
   Choose $\vert k(e)\vert$ distinct
   points on each
edge $e$ of $X$ not lying on $\partial X$.  Provide these points with
positive
coorientation on $e$    if
$k(e)>0$   and with negative coorientation on $e$ if $k(e)<0$ (recall
that $e$
is
oriented). The boundary of each 2-cell of
$X$ meets a certain number   of these  distinguished points. By the
cocycle
condition, their algebraic number   is 0 so that we can join these
points in the
2-cell by
disjoint cooriented intervals  compatible with  the coorientation at
the
endpoints. Proceeding in this way in all 2-cells of $X$
we obtain a cooriented regular graph $\Gamma\subset X$. It is clear
that
$\Gamma$
represents
$s$.   Therefore $\vert \vert
s\vert\vert\leq   \chi_- (\Gamma)=-\chi (\Gamma)=\vert k\vert$.
Conversely,  any cooriented regular graph $\Gamma\subset
X\backslash X^{(0)}$ representing $s$ defines a $\Bbb
Z$-valued   1-cocycle $k$ on $(X, \partial X)$ whose value   on an
(oriented) edge $e$ is equal to the intersection number
$e\cdot \Gamma$. This cocycle represents $s$ and an easy computation
shows that
  $\chi_- (\Gamma)=-\chi (\Gamma)\geq \vert k\vert$. Therefore
  $\vert \vert s\vert\vert\geq  \min_k \vert k\vert$.
 
The formula $\vert \vert s\vert\vert=\min_k \vert k\vert$  is
especially useful
in
the cases where
   either all 0-cells of $X$ lie on $\partial X$ or  $\partial
X=\emptyset$ and
$X$ has only one 0-cell. In both cases every cohomology class $s\in
H^1(X,
\partial X) $ is represented by  a unique cocycle.
 
   \skipaline{\bf 1.7\qua Examples}\qua
   (1)\qua  If $X$ is a compact surface then all elements of  $ H^1(X,
\partial X)$
are
represented by regular graphs consisting of  disjoint embedded circles.
Therefore the
norm $\vert\vert ...\vert\vert_X$ on  $ H^1(X, \partial X;\Bbb R)$
vanishes.

  (2)\qua Let  $\Gamma$ be  a   graph such that all its vertices are incident
to at
least two edges (counted with multiplicity). Let $f$ be a homeomorphism
of
$\Gamma$ onto
itself. The mapping torus, $X$, of $f$ is a 2-complex with void
boundary.
The fibers of the natural fibration $X\to S^1$ determine a class, $s\in
H^1(X)$.
Clearly,  $\vert\vert s\vert\vert \leq \chi_-(\Gamma)$. We   show in
Sect.
3
that
$\vert\vert s\vert\vert= \chi_-(\Gamma)$. This example can be
generalised to
maps
$\Gamma\to \Gamma$ whose mapping torus is a 2-complex.

(3)\qua Let $\Gamma$ be a graph as in the previous example.  The cylinder
$X=\Gamma
\times
[-1,1]$ is
a finite 2-complex with  $\partial X=\Gamma \times \{-1,1\}$. The graph
$\Gamma\times
0 \subset X$ endowed with a coorientation
represents a certain $s\in  H^1(X, \partial X ) $.  The cylinder $X$
embeds in
$\Gamma\times S^1$ in the obvious way and therefore  it follows from
  the   previous example that  $\vert\vert s\vert\vert=
\chi_-(\Gamma)$.

  \skipaline {\bf 1.8\qua Two-complexes associated with group
presentations}\qua Let $\pi$ be a group presented by a finite number  of
generators
and relations
  $\langle x_1,...,x_m : r_1 ,...,r_n \rangle$ where $r_1,...,r_n$ are
words in
the alphabet $  x_1^{\pm 1},...,  x_m^{\pm 1}$. In this subsection we
consider
only
presentations such that each generator appears in the relations at
least twice.
The presentation $\langle x_1,...,x_m : r_1 ,...,r_n \rangle$ gives
rise in the
usual way to a 2-dimensional CW-complex $X$ with one 0-cell, $m$
one-cells and
$n$
two-cells. Let $\#(x_i)$ be the total number of appearances of   $x_i$
in the
words
$r_1,...,r_n$. (A power $x_i^k$ appearing in these words contributes
$\vert
k\vert$ to $\#(x_i)$). It is clear that $\#(x_i)$ is the number of
2-cells of
$X$
adjacent to the $i$-th 1-cell of $X$. By assumption, $\#(x_i)\geq 2$
for all $i$
so
that $\partial X=\emptyset$. Using the formulas of Sect. 1.6, we can
compute
the   norm $ \vert \vert
...  \vert\vert_X $ on $H^1(X;\Bbb R)=H^1(\pi;\Bbb R)$   by
$\vert \vert s\vert\vert_X =\sum_{i=1}^n (\#(x_i)/2-1) \vert
s(x_i)\vert$
for any $s\in H^1(\pi;\Bbb R)$.
This norm    depends on the   presentation of $\pi$. It is   easy to
increase
this norm for instance by adding a tautological relation $x_i
x_i^{-1}=1$.

  We
say that a finite presentation of $\pi$ by generators and relations is
{\sl
minimal}
if the corresponding norm on $H^1(\pi;\Bbb R)$ (considered as a
function) is
smaller than or equal to the  norm on $H^1(\pi;\Bbb R)$ determined by
any other
finite presentation of $\pi$.
For instance, if each generator   appears in the relations exactly
twice, then
the corresponding norm is zero and  the group presentation is minimal.
Another
example:   $\pi =\langle x, y : x^p y^q=1 \rangle$ where $p, q\geq 2$
are
relatively prime integers. A generator   $s\in H^1(\pi)=\Bbb Z$
takes  values $-q$ and $p$ on   $x,y$, respectively.
  The norm of $s$ with respect to this presentation equals  $    (p/2-1)
q+
(q/2-1) p=
  pq-p-q$. We shall show in Sect. 3 that this presentation   is minimal.

  \skipaline {\bf 1.9\qua A related construction}\qua We describe a
related
construction which derives a norm on the cohomology of a compact
surface
$\Sigma$ from
a  family  of  loops on $\Sigma$.
  Let  $\alpha=\{\alpha_i\}_i$ be  a finite family  of closed curves in
$\Int
\Sigma$
whose all
crossings and  self-crossings   are
transversal double intersections. We  define  a norm $\vert
...\vert_\alpha$  on
$
H^1(\Sigma,
\partial \Sigma;\Bbb R)$ as follows. For  any $s\in
   H^1(\Sigma, \partial \Sigma ) $, set
$\vert s\vert_\alpha=\min_S \# (S\cap \cup_i \alpha_i)$ where $S$ runs
over
cooriented
closed
1-dimensional submanifolds  of $\Sigma $ representing $s$ and meeting $
\cup_i
\alpha_i$ transversely
   (in the complement of the set of double points of $ \cup_i
\alpha_i$). Here $ \# (S\cap \cup_i \alpha_i)$ is
the number of points in  $S\cap \cup_i \alpha_i$.
It is easy to check that  $\vert ...\vert_\alpha$ is a   homogeneous
norm on
$ H^1(\Sigma, \partial \Sigma )$.
As usual it extends uniquely  to a  homogeneous continuous  norm, also
denoted
$\vert
...\vert_\alpha$,  on  $ H^1(\Sigma, \partial \Sigma;\Bbb R)$. This
norm is
preserved
under the
first and third Reidemeister moves on the loops $\{\alpha_i\}_i$ but in
general
is not
preserved under the second Reidemeister move.
  A simple
example is provided by a small loop $\alpha\subset \Sigma$ bounding a
disc in
$\Sigma$.
The   norm
$\vert ...\vert_\alpha$  on $ H^1(\Sigma, \partial \Sigma;\Bbb R)$ is
zero. On
the
other
hand
we can deform $\alpha$ into an immersed loop $\beta$ in  $S$ which
splits
$\Sigma$
into 2-discs. The   norm
$\vert ...\vert_\beta$ is then non-degenerate.

  The norm $\vert ...\vert_\alpha$  on $ H^1(\Sigma, \partial
\Sigma;\Bbb R)$ is
related
to the
norm on the 1-cohomology of 2-complexes   as follows.
Let
   $\alpha=\{\alpha_i\}_i$ be  a finite family  of loops  in $\Int
\Sigma$ as
above.
Let $X$ be the 2-complex obtained by gluing 2-discs to $\Sigma$ along
these
loops. It is
clear that $\partial X=\partial \Sigma$. We can identify  $ H^1(X,
\partial
X;\Bbb R)$
with the linear subspace of  $ H^1(\Sigma, \partial \Sigma;\Bbb R)$
consisting
of
cohomology classes whose evaluation on  the  loops  $\{\alpha_i\}_i$ is
$0$. Then
the norm $\vert \vert ...\vert\vert_X$
on  $ H^1(X, \partial X;\Bbb R)$ is the restriction of
$(1/2)\vert...\vert_\alpha$. Indeed,  any  regular
graph
$\Gamma\subset X$ consists of a   closed    1-manifold
$S=\Gamma\cap
\Sigma  $
and several intervals lying in the glued 2-discs and connecting the
   points of
    $S\cap \cup_i \alpha_i$.
  All vertices of $\Gamma$ are trivalent and therefore
$\chi_-(\Gamma)=-\chi(\Gamma)=
(1/2)\, \# (S\cap \cup_i \alpha_i)$.
The 2-complexes   obtained  in this way from homotopic systems of loops
are
  (simply) homotopy equivalent. The example  above implies that  the
norm
  $ \vert \vert ...\vert\vert$ in general is not preserved under
(simple)
  homotopy equivalences of 2-complexes.

\section{The Alexander-Fox
polynomials and norms}

The Alexander polynomial is mostly known in the context of knot
theory.  Fox observed that this polynomial depends only on the knot
group and in fact can be defined for an arbitrary finitely generated
group. In this section we recall the relevant definitions following
\cite {Fox}. In generalisation of the standard Alexander-Fox
polynomial, we introduce twisted Alexander-Fox polynomials and
consider the associated norms on 1-cohomology of groups.
 
  Fix throughout this section a finitely generated group $\pi$.  Set
$H=H_1(\pi)$
and $
G=H/\Tors
H$. The ring homomorphism     $\Bbb  Z[H ]\to \Bbb  Z[G ]$ induced by
the
projection $H\to G$ will be denoted by
$\proj$.

\skipaline {\bf 2.1\qua  The   elementary
ideals}\qua
  The group $\pi$ determines an
increasing
sequence of ideals  $E_0(\pi)\subset E_1(\pi)\subset E_2(\pi)\subset
...  $ of
the
group ring $ \Bbb  Z[H]$
called {\it the elementary
ideals} of $\pi$. They can be computed from an arbitrary presentation
of $\pi$
by
generators and relations
$\langle x_1,...,x_m: r_1,r_2,...\rangle$ with finite $m\geq 1$.  Here
each
$r_i$ is viewed as an
element of the free group,
$F$, generated by $x_1,...,x_m$; the number of relations can be
infinite.
Every
$f\in F$ can be uniquely expanded in $\Bbb  Z[F]$ as $1+\sum_{j=1}^m
f_{j}
(x_j-1)$
with $f_{1},..., f_{m}\in \Bbb  Z[F]$.  The element $f_{j}\in \Bbb
Z[F]$ is
called
the {\sl $j$-th Fox derivative} of $f$ and denoted by $\partial f/
\partial
x_j$.
Consider the matrix $ [\partial r_i/ \partial x_j]_{i,j}$ over $\Bbb
Z[F]$.
Applying the natural projections $\Bbb  Z[F] \to \Bbb  Z[\pi]\to \Bbb
Z[H]$
to the entries of this matrix we obtain a matrix, $A$, over $\Bbb
Z[H]$  called the {\it Alexander-Fox matrix} of the presentation
$\langle x_1,...,x_m: r_1,r_2,...\rangle$.  It has $m$ columns and
possibly
infinite
number
of
rows.  Adding if necessary to   $r_1,r_2,...$ several copies of the
neutral element $1\in F$ we can assume that $A$ has at least $m$ rows.
For
$d=0,1,...$, the ideal $E_d(\pi)\subset \Bbb  Z[H]$ is generated
by the
minor determinants of $A$ of order $m-d$.  This ideal does not depend
on the
presentation of $\pi$.  We shall be interested only in the ideal
$E_1(\pi)$
which
will be denoted $E(\pi)$.

\skipaline {\bf 2.2\qua The Alexander-Fox polynomials}\qua
Consider the ideal $ \proj (E(\pi))\subset \Bbb  Z[G ]$.  Since $\Bbb
Z[G ]$
is
a
unique factorization domain, one can consider the greatest common
divisor of
the
elements of $\proj (E(\pi))$.  This $\gcd$  is an element of $ \Bbb
Z[G ]$
defined
up to multiplication by $\pm G $.  It is called the {\sl Alexander-Fox
polynomial}
of
$\pi $ and denoted $\Delta(\pi)$.

    The obvious inclusion $\proj (E(\pi))\subset  \Delta(\pi)\, \Bbb
Z[G]$ can be
slightly
improved provided $\rk\, H\geq 2$.  Namely, if $\rk\, H\geq 2$, then
$$\proj
(E(\pi))\subset  \Delta(\pi)\,  J\eqno (2.a)$$ where $J$ is the
  augmentation
ideal of   $\Bbb Z[G]$. This inclusion goes back to   \cite {Fox}, Prop. 6.4
at
least in the case $\Tors H=0$. We give a proof of (2.a) at the end of
Sect. 2.

   In generalisation of $\Delta(\pi)$, we define {\sl twisted
Alexander-Fox
polynomials}
of $\pi$ numerated by $\sigma\in (\Tors H)^*=\Hom (\Tors H, \C^*)$.
Fix a
splitting $H=\Tors H \times G$.  For $\sigma\in (\Tors H)^*$, consider
the ring
homomorphism $\tilde \sigma:\Bbb  Z[H] \to \C [G]$ sending $fg$
with $f\in
\Tors H, g\in G$ to $\sigma(f) g$ where $\sigma(f)\in \C^*\subset
\C$.
The ring $\C  [G ]$ is a unique factorization domain and we can set
$\Delta^\sigma (\pi)= \gcd  \, \tilde \sigma (E(\pi))$.  This $\gcd$
is an
element of $ \C [G]$ defined up to multiplication by elements of
$G$ and
nonzero complex numbers.  Under a different choice of the splitting
$H=\Tors H
\times G$, the polynomial $\Delta^\sigma(\pi)$, represented say by
$\sum_{g\in
G} c_g
g$ with
$c_g\in
\C$, is replaced by $\sum_{g\in G} c_g \sigma (\psi (g)) g$ where
$\psi\in
\Hom
(G, \Tors H)$.  For $\sigma=1$, we have $ \Delta^1 (\pi)=\C^*
\Delta (\pi)
$.

\skipaline {\bf 2.3\qua The Alexander-Fox polytopes and
norms}\qua   Fix
$\sigma\in
(\Tors H )^* $. In analogy with
the Newton polytope of a polynomial, we can derive from
$\Delta^\sigma(\pi)$  an
{\sl
Alexander-Fox polytope}  (or briefly AF-polytope) $P^\sigma(\pi)\subset
H_1(\pi;\Bbb  R)$.  Pick   a
representative
$\sum_{g\in G } c_g g \in \C[G]$ of  $\Delta^\sigma(\pi)$.
Set $$ P^\sigma(\pi) = {HULL} \, (\{\,\frac {1}{2}  ({g^{real}
-(g')^{real}} )
\,\,\vert \,\,
g,g' \in G
,\, \, \, c_g\neq 0, \, c_{g'} \neq 0\,\})    $$
where $g^{real}\in H_1(\pi;\Bbb  R)$ is the real homology class
represented by
$g\in G$ and for a subset $S  $ of a linear space,   $HULL(S)$ denotes
the
convex
hull of
$S$.
By convention, if
$\Delta^\sigma(\pi)=0$ then $P^\sigma(\pi)=\{0\}$.
The polytope $P^\sigma(\pi)$ is a compact convex
polytope
symmetric
in the origin and independent of     the representative
$\sum_{g } c_g g$.      Its vertices lie on the half-integral lattice
$(1/2) G$
where $G\subset H_1(\pi;\Bbb  R)$ consists of integral homology
classes.

We define the {\sl
Alexander-Fox norm}  (or briefly AF-norm) $\vert \vert...\vert
\vert^\sigma$ on
$H^1(\pi;\Bbb
R)$ by
$$\vert \vert s\vert \vert^\sigma=  2 \max_{x\in P^\sigma(\pi)} \vert
s(x) \vert =
\max_{g,g'\in G, c_g  c_{g'} \neq 0}
\vert
s(g) - s(g') \vert  $$
where $s\in H^1(\pi;\Bbb  R)$  and $s(x)\in \Bbb R$ is the
evaluation
of $s$ on $x$.  This
  norm   is   continuous and homogeneous. It
was
first considered  in the case $\sigma=1$  by C. McMullen \cite {McM}.

The AF-norms are natural with respect to group isomorphisms:
For a
group isomorphism $\varphi:\pi'\to \pi$ and   $s\in H^1(\pi;\Bbb R),
\sigma\in
(\Tors H_1(\pi) )^*$, we have $\vert \vert s \vert \vert^\sigma
=\vert \vert \varphi^*(s) \vert \vert^{\sigma \varphi_*}$
where $\varphi^*$ and $\varphi_*$ are  the induced homomorphisms
$H^1(\pi;\Bbb R)
\to H^1(\pi';\Bbb R)$ and $\Tors H_1(\pi') \to \Tors H_1(\pi) $,
respectively.

   \skipaline {\bf 2.4\qua Examples}\qua
   (1)\qua If    $\pi$ has a presentation with   $m$
generators
and
$\leq m-2$ relations then  $E(\pi)=0$
and
the AF-norms
on $H^1(\pi;\Bbb  R)$ are   0.
 
  (2)\qua  Let $p, q\geq 2$ be relatively prime integers and $\pi=\langle
x, y: x^p y^q=1\rangle$. Let $t$ be a generator of
$H_1(\pi)=\Bbb Z$. Set $n=pq-p-q+1$. The  polynomial   $\Delta(\pi)$ is
represented by the
Laurent polynomial
  $(t^{pq}-1) (t-1) (t^p-1)^{-1} (t^q-1)^{-1}$ with lowest term $1$ and
highest
term
  $t^{n}$.
  The AF-polytope in $H_1(\pi;\Bbb  R)=\Bbb R$ is the interval with
endpoints $
  -(n/2) t^{real}$ and $ (n/2)  t^{real}$.
  The AF-norm of both generators of $H^1(\pi)=\Bbb Z$ is equal to
  $n=pq-p-q+1$.
 
  (3)\qua Let
  $\pi=\langle x,y :
x^k y^l x^{-k} y^{-l}=1, y^m=1\rangle$
where  $k,l\geq 1, m\geq 2 $.
   It is clear that
$H_1(\pi)= \Bbb Z \times (\Bbb Z/m \Bbb Z)$ with generators $[x], [y]$
represented
by $x,y$.
A direct computation  shows   $E(\pi)$ is generated
by 3 elements: $1+[y]+...+ [y]^{m-1}, (1+[x]+...+ [x]^{k-1})
([y]^l-1)$, and
$([x]^k-1)
(1+[y]+...+ [y]^{l-1})$. Setting $[y]=1$ we obtain  that
$\Delta(\pi)=\gcd (l,m)$.
The corresponding AF-norm is zero.
  Let $\sigma$  be a nontrivial character   of $\Tors H_1(\pi)=\Bbb Z/m
\Bbb Z$.
  Then $\zeta=\sigma([y])\neq 1$
  is a complex root of unity   of order $m$.
  If $\zeta^l=1$ then $\Delta^{\sigma } (\pi)=0$ and the corresponding
AF-norm is
zero. If $\zeta^l\neq 1$ then $\Delta^{\sigma } (\pi)=1+t +...+t^{k-1}$
where
$t$ is the generator of $H_1(\pi)/\Tors H_1(\pi)$ represented by $x$.
The
corresponding AF-norm of both generators of
$H^1(\pi)=\Bbb Z$ is equal to
  $k-1$.  This example shows that the twisted AF-polynomials may provide
more
interesting norms than the untwisted AF-polynomial.

\skipaline {\bf 2.5\qua Remark}\qua   The structure of the   ideal
$E(\pi)$ can be sometimes described using the theory of   Reidemeister
torsions.
Suppose that $\pi=\pi_1(X)$ where $X$ is
a finite connected 2-complex   with $\chi (X)=0$.
As above, set $  H=H_1(X)=H_1(\pi), G=H/\Tors H$.
The  maximal abelian torsion  $\tau(X)$ is an element of the
commutative ring $Q(H)$ obtained from   $\Bbb Z[H]$ by inverting all
non-zerodivisors (see \cite {Tu1},\cite {Tu2}). The natural homomorphism $\Bbb
Z[H]\to
Q(H)$
is an
inclusion and we can identify $\Bbb Z[H]$ with its image. Then $E(\pi)=
\tau(X) I$
where $I$ is the augmentation ideal of $\Bbb Z[H]$ (for a proof, see
\cite {Tu1}, p.
689). If $\rk H\geq 2$, then $\tau(X)\in \Bbb Z[H]$ and for any
$\sigma\in (\Tors H)^*$, the twisted AF-polynomial $\Delta^\sigma(\pi)$
is
represented by
$\tilde \sigma(\tau (X))\in \C [G]$. If $\rk H=1$, then
$\tau(X)$ splits as  a sum $a  +  (t-1)^{-1} \Sigma$ where $a\in \Bbb
Z[H]$,
$\Sigma=\sum_{f\in \Tors H} f \in \Bbb Z[H]$, and $t$ is any element of
$H$
whose
projection $\proj (t)\in G=\Bbb Z$ is a generator. Then    for any
non-trivial
character  $\sigma\in
(\Tors H )^*$, the  polynomial $\Delta^\sigma(\pi)$ is
represented by $\tilde \sigma(a)\in \C [G]$. The    polynomial
$\Delta (\pi)$ corresponding to $\sigma=1$ is represented by $\proj
((t-1)a) +
\vert \Tors H \vert\in \Bbb Z [G]$.

\skipaline {\bf 2.6\qua   Proof of (2.a)}\qua Consider a
presentation
$\langle x_1,...,x_m: r_1,r_2,...\rangle $ of
$\pi$ by
generators and relations
  with finite $m\geq 1$ and at least $m$
relations.
  Let $A$ be the Alexander matrix of this presentation. It is enough to
show
that for
any
  minor determinant  $D$ of $A$ of order $m-1$,
we have $\proj (D) \in \Delta(\pi) J$. Assume for concreteness that $D$
is the
determinant
of a
submatrix of the first $m-1$ rows of $A$.
Let $\pi'$ be the group $\langle x_1,...,x_m:
r_1,r_2,...,r_{m-1}\rangle$. Set
$H'=H_1(\pi')$. The
natural
surjection $
H'\to H=H_1(\pi)$   induces a ring
homomorphism
$\Bbb Z[H']\to \Bbb Z[H ]$  denoted   $\psi$.
It follows from definitions  that $D\in \psi  (E(\pi')) \subset
E(\pi)$.
Note that $ \rk H'\geq
\rk H\geq
2$.

Consider the
2-dimensional CW-complex $X$
determined
by the   presentation $\langle x_1,...,x_m: r_1,r_2,...,
r_{m-1}\rangle$.
Clearly, $ \pi_1(X)=\pi'$  and $\chi (X)=0$. By Remark 2.5,
$E(\pi')= \tau I'$ where  $\tau \in \Bbb Z[H']$ and  $I'$ is
the augmentation ideal of
$\Bbb
Z[H']$.
Applying $ \proj \circ \psi $
we obtain that
$$(\proj \circ \psi)  (\tau) J = (\proj \circ \psi)  (\tau I')= (\proj
\circ \psi)
  (E(\pi'))   \subset \proj
(E(\pi))
\subset \Delta(\pi)\, \Bbb Z[G].$$
Since $\rk H\geq
2$, we have $\gcd J=1$ and hence
$\Delta(\pi)$ is a divisor of $(\proj \circ \psi) (\tau)\in \Bbb Z[G]$.
Therefore
  $\proj (D)\in (\proj \circ \psi)
(E(\pi'))=
(\proj \circ \psi)  (\tau) J \subset
\Delta(\pi) J$.

\section{Main theorem}

To state our main theorem  it is convenient to introduce a   {\sl
trivial norm}
$\vert
...\vert_0
$ on the real 1-cohomology
$H^1(X;\Bbb  R)$ of any CW-space $X$. If the first Betti number of $X$
is $ \neq
1$
then $\vert s\vert_0=0$ for all $s\in H^1(X;\Bbb  R)$.
If the first Betti number of $X$ is $1$, then
  $\vert...\vert_0$ is the unique homogeneous norm
  on $H^1(X;\Bbb  R)$ taking value $ 1$ on both generators of
$\Bbb  Z=H^1(X;\Bbb  Z) \subset H^1(X;\Bbb  R) $.

  \skipaline {\bf 3.1\qua Theorem}\qua {\sl  Let $X$ be a
connected finite
2-complex
with
$\partial X=\emptyset$.
  For any $s\in H^1(X;\Bbb  R)$ and any $\sigma \in (\Tors H_1(X) )^*$,
$$\vert\vert s
\vert\vert_X\geq  \vert \vert s\vert
\vert^\sigma -
\delta_\sigma^1 \vert s \vert_0   \eqno (3.a)$$ where $\vert \vert
...\vert
\vert^\sigma$ is
the Alexander-Fox  norm  on $H^1(X;\Bbb  R)=H^1(\pi_1(X);\Bbb  R)$
determined
by
$\sigma$ and
$\delta_\sigma^1=1$ if $\sigma=1$ and
$\delta_\sigma^1=0$ otherwise. }

\skipaline Theorem 3.1 will be proven in Sect. 4. Note that the norm
$\vert \vert
... \vert
\vert^\sigma$ on $H^1(X;\Bbb  R)$ does not depend on the choice of a
base point in
$X$ because of the invariance of the AF-norms under group isomorphisms.
In the
case $\rk
H_1(X)\geq 2$,   (3.a) simplifies to
$\vert\vert s
\vert\vert_X\geq    \vert \vert s\vert
\vert^\sigma  $.

Inequality (3.a) has a version for 1-cohomology classes on
  3-manifolds, where on
the left hand side  appears the Thurston norm of $s$ and
   the right hand side
is $\vert \vert s\vert
\vert^\sigma -
2\delta_\sigma^1 \vert s \vert_0$. The author plans to discuss this
version of
Theorem 3.1  elsewhere.

   \skipaline {\bf 3.2\qua Corollary}\qua {\sl Let $M$ be a
connected
manifold (possibly with boundary)  of dimension $ \geq 3$. Let $X$ be a
connected
finite
2-complex
with
$\partial X=\emptyset$   embedded in $M$ such that the inclusion
homomorphism
$\pi_1(X)\to \pi_1(M)$ is an isomorphism. Let  $S\subset M$ be a
cooriented
compact submanifold of $M$ of codimension 1 intersecting $\partial M$
along
$\partial S$ and intersecting $X$
transversely along a
regular graph $\Gamma=S\cap X$. Let $s\in H^1(M;\Bbb Z)$ be the
cohomology class
represented by
$S$.
Then  $$\vert \chi (\Gamma) \vert\geq \max_{\sigma \in (\Tors H_1(M)
)^*}
(\vert
\vert s\vert \vert^\sigma -
\delta_\sigma^1 \vert s \vert_0 ) \eqno (3.b)$$ where $\vert \vert
...\vert
\vert^\sigma$ is
the Alexander-Fox norm on $H^1(M;\Bbb  R)$ determined by $\sigma$.}

\skipaline  The assumption $\pi_1(X)= \pi_1(M)$
ensures   that $\pi_1(M)$ is finitely generated so that the AF-norms on
$H^1(M;\Bbb  R)$ are well defined.

  To deduce Corollary 3.2 from Theorem 3.1,  set $s'= s\vert_X \in
H^1(X;\Bbb
R)$. Clearly, $s'=s_\Gamma$.
Therefore   $ \vert \chi (\Gamma) \vert = \chi_- (\Gamma)\geq \vert
\vert
s'  \vert \vert$. By Theorem 3.1 and the assumption $\pi_1(X)=
\pi_1(M)$,
$$\vert
\vert s' \vert \vert
\geq \max_{\sigma \in (\Tors H_1(X) )^*} (\vert \vert s'\vert
\vert^\sigma -
\delta_\sigma^1 \vert s' \vert_0 )=\max_{\sigma \in (\Tors H_1(M) )^*}
(\vert
\vert
s\vert \vert^\sigma -
\delta_\sigma^1 \vert s \vert_0 ).$$

The statement of Corollary 3.2 is not   specific  about the category
of
manifolds. In fact the corollary  extends to a much broader setting
where $M$
  is an arcwise connected
space and $S$ is a subspace of  $M$
which has a cylinder
neighborhood   $U=S\times [-1,1]\subset M$ such that $S=S\times 0$.
A coorientation
of $S$
is defined in the obvious way and
determines
a  cohomology class $s\in H^1(M)$ as in Sect. 1.2.
As above, $X\subset M$ is a
connected
finite
2-complex
with
$\partial X=\emptyset$     such that the inclusion
homomorphism
$\pi_1(X)\to \pi_1(M)$ is an isomorphism and
$U\cap
X=\Gamma
\times [-1,1]$ where $\Gamma=S\cap X$ is a graph in $X$.  Then   we
have     (3.b).

Corollary 3.2 can be applied in various geometric situations. For
instance, if
$M$ is a compact triangulated manifold of dimension $\geq 3$ then we
can take
$X$
to be the 2-skeleton
  of $M$. If $M$ is a compact 3-manifold then we can take   $X$ to be a
spine of
$M$
or a spine of punctured $M$.

   \skipaline {\bf 3.3\qua Corollary}\qua {\sl Let $\pi$ be a
group presented
by a finite number  of generators
and relations
  $\langle x_1,...,x_m : r_1 ,...,r_n \rangle$ where $r_1,...,r_n$ are
words in
the alphabet $  x_1^{\pm 1},...,  x_m^{\pm 1}$ such that (in the
notation of Sect.
1.8) $\#(x_i)\geq 2$ for $i=1,...,m$. Then for any $s\in H^1(\pi;\Bbb
R)$,
$$\sum_{i=1}^n \, (\#(x_i)/2-1)\,  \vert s(x_i)\vert\geq \max_{\sigma
\in (\Tors
H_1(\pi
) )^*} (\vert \vert s\vert
\vert^\sigma -
\delta_\sigma^1 \vert s \vert_0 ). \eqno (3.c)$$
}

\skipaline This corollary is  obtained by an application of Theorem 3.1
to the
2-complex determined by the presentation $\langle x_1,...,x_m : r_1
,...,r_n
\rangle$.

  \skipaline {\bf 3.4\qua Examples}\qua (1)\qua The computations  in
Sect. 1.8 and
2.4.2 show that for
  the group presentation
$\langle
x, y: x^p y^q=1\rangle$ the inequality (3.c) is an equality. Thus, this
presentation is minimal in the sense of Sect. 1.8.
(It would be interesting to extend this fact   to groups of other
fibered knots
in
$S^3$).

(2)\qua    Let
  $\pi=\langle x,y :
x^k y^l x^{-k} y^{-l}=1, y^m=1\rangle$
where  $k,l\geq 1, m\geq 2 $. We claim that if $m$ does not divide $l$
then this
presentation is minimal. Indeed, there is a nontrivial character
$\sigma$  of
$\Tors H_1(\pi)=\Bbb Z/m \Bbb Z$ such that $(\sigma([y]))^l\neq 1$.
  The computations in Sect. 2.4.3 show that $\vert \vert
s\vert\vert^\sigma=k-1$
where $s$ is a generator of $H^1(\pi) =\Bbb Z$.
  The left hand side of (3.c) is
  $(2k)/2-1=k-1$. Hence (3.c) is  an equality for this  presentation of
$\pi$ which is
therefore minimal. Examples 1 and 2 show   that the
estimate in Theorem 3.1 is sharp.

(3)\qua Consider $\Gamma, f, X, s$ from the mapping
   torus of Example 1.7.2.  We  will  deduce the equality $\vert\vert
s\vert\vert= \chi_-(\Gamma)$ from  Theorem
3.1.
We need only to prove that $\vert\vert s\vert\vert\geq
\chi_-(\Gamma)$.
It is enough to consider the case of connected $\Gamma$.
By (3.a), it is enough to show that $\chi_-(\Gamma) =\vert \vert s\vert
\vert^1 -
   \vert s \vert_0$.  To this end we shall compute the (untwisted)
AF-polynomial
$\Delta(X)=\Delta(\pi_1(X))$.  We can  deform $f:\Gamma\to
\Gamma$ so that it fixes a point  $\gamma\in \Gamma$. Let $x_1,...,x_n$
be  free
generators of the free group $\pi_1(\Gamma, \gamma)$ where $n\geq 1$.
The group
$\pi_1(X)$ can
be
presented by $n+1$ generators $x_1,...,x_n, T$ subject to $n$ relations
$ T x_i
T^{-1} (f_{\#} (x_i))^{-1}=1$ where $i=1,...,n$ and $f_{\#}$ is the
endomorphism
of
$\pi_1(\Gamma, \gamma)$ induced by $f$.
   Set $G=H_1(X)/\Tors H_1(X)$ and let $G'$ be the corank 1 sublattice
of $G$
generated by the classes $[x_1],..., [x_n]\in G$ of $x_1,...,x_n$. Set
$\delta=1$
if $\rk G\geq 2$
   and  $\delta=0$ if $\rk G=1$.
  A direct computation using the Fox differential calculus   gives
   $\Delta(X)= (t-1)^{-\delta} \det (tE_n -A)$ where $t=[T]\in G$ is the
class of
$T$, $E_n$ is the
unit
$(n\times n)$-matrix, and $A$ is the $(n\times n)$-matrix over   $\Bbb
Z[G']$
obtained as the image of the matrix $(\partial f_{\#} (x_i)/\partial
x_j)_{i,j=1,...,n}$ under the natural ring homomorphism $\Bbb
Z[\pi_1(\Gamma,
\gamma)]\to \Bbb Z [G']$ sending each $x_i$ to   $[x_i]$. Clearly,
$\det (tE_n
-A)= a_0+ a_1t+...+a_{n-1} t^{n-1} +t^n$ where $a_0,..., a_{n-1}\in
\Bbb Z[G']$.
   Since $f_{\#}$ is an isomorphism, the sum of coefficients of $a_0=\pm
\det A \in \Bbb Z[G'] $ is $\pm 1$   and therefore $a_0 \neq 0$.
   By definition, $s(G')=0$ and $s(t)=\pm 1$. If   $\rk G\geq 2$ then
$\Delta(X)=  -a_0+  ...+
t^{n-1}
$ and   $\vert \vert s\vert
\vert^1 -
   \vert s \vert_0=\vert \vert s\vert
\vert^1=n-1=\chi_-(\Gamma)$. If $\rk G=1$ then $\Delta(X)=  a_0+  ...+
t^{n}   $ and
$\vert \vert s\vert
\vert^1 -
   \vert s \vert_0=n-1=\chi_-(\Gamma)$.

\section{Proof of Theorem 3.1}

{\bf 4.1\qua Preliminaries on modules}\qua  Let
$\Lambda$ be a
commutative ring with unit.  For a finitely generated $\Lambda$-module
$X$
consider
a $\Lambda$-linear homomorphism  $ f:\Lambda^n \to \Lambda^m $
with finite
$m$ and  $\Coker f= X$.  The {\it $i$-th elementary ideal}
$E_i(X)\subset \Lambda$ with
$i=0,1,...$
is generated by
the
$(m-i)$-minors of the matrix of $f$. If $m-i>n$ then $E_i(X)=0$; if
$m-i\leq 0  $,
then $E_i(X)=\Lambda$.  The ideal
$E_i(X)$ is independent of $f$.
If $\Lambda$ is a unique factorization domain then the {\it $i$-th
Alexander
invariant} $\Delta_i (X)\in \Lambda$ of $X$ is the greatest common
divisor of
the
elements of $E_i(X)$.  It is well-defined up to multiplication by units
of
$\Lambda$.     Note that
$\Delta_i(X\oplus \Lambda)=\Delta_{i-1} (X)$ for $i\geq 1$.

   Let $\Lambda=\C [t^{\pm 1}]$.  The Alexander invariants of a
finitely
generated $\Lambda$-module $X$ can be computed as follows.  Since
$\Lambda$ is
a
principal ideal domain, $X=
\oplus_{r=1}^m
(\Lambda/\lambda_r)$ where $\lambda_1,..., \lambda_m\in \Lambda$ and
$\lambda_{i+1}$ divides $\lambda_i$ for all $i$.  Then $
\Delta_i(X)=\prod_{r=i+1}^m\lambda_r$ for $i<m$ and $\Delta_i(X)=1$ for
$i\geq
m$.
  The maximal $r$ such that
$\lambda_1=...=\lambda_r=0$ is called the {\sl rank} of $X$ and denoted
$\rk_\Lambda\, X$.
It
is
clear that $\rk_\Lambda\, X=\dim_{\Bbb  Q(t)}\,(\Bbb  Q(t) \otimes
_\Lambda X)$
where $\Bbb  Q(t)$
is
the field of fractions of $\Lambda$.
We have $\rk_\Lambda\, X=0 \Leftrightarrow \Bbb  Q(t) \otimes _\Lambda
X= 0
\Leftrightarrow
\Delta_0(X)\neq 0 $.  If $\rk_\Lambda\, X=0$ then $X$  is a finite
dimensional
$\C$-linear space and $\dim_{\C}\, X=
\rm
{span}
\Delta_0(X)  $ where the span of a nonzero Laurent polynomial $
\sum_n
a_n
t^n\in \Lambda$ is $ \max_{m,n, a_ma_n\neq 0} \vert m-n\vert$.

   \skipaline {\bf 4.2\qua Preliminaries on twisted
homology}\qua
We recall the notion of twisted homology.  Let $X$ be a
connected
CW-space and $H=H_1(X)$.  Let  $\Lambda$ be  a commutative ring with
unit
and
$\varphi$
be  a ring homomorphism $\Bbb Z[H] \to \Lambda$.  We view $\Lambda$ as
a (right)
$\Bbb
Z[H]$-module via $ \lambda z =  \lambda \varphi(z)$ for $
\lambda
\in
\Lambda, z \in \Bbb Z[H]$. Let $p:\hat X\to X$ be  the maximal abelian
covering
of
$X$ (with induced CW-structure) corresponding to the commutant of
$\pi_1(X)$. The
action
of $H$ on $ \hat X $ by
deck transformations makes  the cellular chain complex
$C_*(\hat X)$ a complex of (free) left $\Bbb Z[H]$-modules. By
definition,
$$H_*^{
\varphi}(X)=H_*( \Lambda\otimes_{\Bbb
Z[H]}
C_*(\hat X)). $$   Note that $H_*^{ \varphi}(X)$ is a
$\Lambda$-module.
The
twisted homology extends to cellular pairs $  Y\subset X$ by
$$H_* ^{
\varphi}(X,Y)=H_*( \Lambda\otimes_{\Bbb Z[H]} C_*(\hat X)/C_*(p^{-1}
(Y)))
$$
where $C_*(p^{-1} (Y))$ is the chain subcomplex of  $C_*(\hat X)$
generated by
cells
of $\hat X$ lying in $ p^{-1} (Y)
$.

The twisted homology is invariant under cellular subdivisions and forms
the
usual
exact homology sequences such as the Mayer-Vietoris homology sequence
and the
homology sequence of a pair.  Using a CW-decomposition of $X$ with one
$0$-cell,
one can check that $H_0^{ \varphi}(X)=\Lambda/ \varphi (I ) \Lambda$
where $I$
is
the augmentation ideal of $\Bbb Z[H]$.

\skipaline {\bf 4.3\qua Preliminaries on weighted graphs}\qua
     The notion of weighted
graphs    formalizes graphs with
parallel
components.  A {\sl weighted graph} in a 2-complex
$X$ is a cooriented regular graph
$\Gamma\subset X$ such that each its component $\Gamma_i$ is endowed
with a
positive integer $w_i$ called the {\sl weight} of $\Gamma_i$. We  write
$\Gamma=\cup_{i  }
(\Gamma_i, w_i)$.
  A weighted graph $\Gamma=\cup_{i  }
(\Gamma_i, w_i)$ in $X $ gives rise to an
(unweighted)  cooriented regular
graph $\Gamma^u\subset X$ obtained by replacing each   $\Gamma_i$
by
$w_i$ parallel copies in a small neighborhood of $\Gamma_i$.  We say
that
$\Gamma$
represents the  cohomology  class
$s_\Gamma=s_{\Gamma^u}=\sum_{i }  w_i s_{\Gamma_i} \in H^1(X) $. Set
$\chi_-
(\Gamma)=\chi_-
(\Gamma^u)=\sum_{i  }  w_i\, \chi_- (\Gamma_i) $.

\skipaline {\bf 4.4\qua Lemma}\qua {\sl  Let $X$ be a
connected finite
2-complex
with
$\partial X=\emptyset$. Every  nonzero $s\in
H^1(X)
$  can be represented by a
weighted
graph $\Gamma\subset X$ such that $\chi_-
(\Gamma)=\vert\vert s \vert\vert$ and   $X\backslash \Gamma$ is
connected.}

\proof
Consider first an arbitrary  weighted graph
$\Gamma=\cup_{i
}
(\Gamma_i, w_i)$
in $X$.  By \lq\lq decreasing the weight of   $\Gamma_i$ by 1" we mean
the transformation which reduces $w_i$ by $1$ and keeps the other
weights.  If
$w_i=1$, then this transformation removes $\Gamma_i$ from $\Gamma$.

Assume that $X\backslash \Gamma$ is not connected. For a  component $N
$ of
$X\backslash \Gamma$, we define a {\it reduction} of
$\Gamma$ along   $N $. Let $\alpha_+$ (resp. $\alpha_-$)  be the set of
all
$i$ such that $N$ is adjacent to $\Gamma_i$ only on the positive (resp.
negative)
side.
  The sets $\alpha_+, \alpha_-$ are disjoint.
Since $N\neq X\backslash \Gamma$,   at least one of
these two sets
   is non-void.
Counting the number of entries and exits in $N$ of a loop on $X$  we
observe that
$\sum_{i\in \alpha_{+}} s_{\Gamma_i} = \sum_{i\in \alpha_{-}}
s_{\Gamma_i}\in
H^1(X)$.
We modify $\Gamma$ as follows.
    If $\alpha_{+}\neq \emptyset$ and $\sum_{i\in \alpha_{+}} \chi_-
(\Gamma_i)\geq  \sum_{i\in \alpha_{-}} \chi_-
(\Gamma_i)$,
then we decrease by 1 the weights of all $\{\Gamma_i\}_{i\in
\alpha_{+}}$ and
increase by 1
the weights of   all $\{\Gamma_i\}_{i\in \alpha_{-}}$.  If $\alpha_{+}=
\emptyset$ or  $\sum_{i\in
\alpha_{+}} \chi_- (\Gamma_i)< \sum_{i\in \alpha_{-}} \chi_-
(\Gamma_i)$
  then we increase by 1 the weights of all $\{\Gamma_i\}_{i\in
\alpha_{+}}$ and
decrease by 1
the weights of   all $\{\Gamma_i\}_{i\in \alpha_{-}}$.    This
   yields another weighted graph $\Gamma'$ such that $s_{\Gamma'}
=s_\Gamma
$ and $\chi_- (\Gamma')\leq \chi_- (\Gamma)$.
    Iterating
this transformation, we eventually remove from $\Gamma$ at least one
component
incident to $ N$ on one side.  Let us call this iteration the
{\sl reduction}
of $\Gamma$ along $N$.  The   reduction does not increase $\chi_- $,
preserves $s_\Gamma$ and
strictly decreases the number of
components of $X\backslash \Gamma$.
  If $\partial \overline N$ is connected then the   reduction along $N$
removes
   $\partial \overline N$ from $\Gamma$.

To   prove the lemma,  represent $s $ by a cooriented regular
graph   $S\subset X$ such that   $\chi_-
(S)=\vert\vert s \vert\vert$.  We   view $S$ as a weighted graph
with  weights of all components equal to $1$. If $X\backslash S$ is
connected
then
$S$ satisfies the requirements of the lemma.  If $X\backslash S$ is not
connected
then
iteratively applying to $S$
reductions
along
   components of $X\backslash S$ we eventually obtain a weighted
graph,
$\Gamma$,
such that $X\backslash \Gamma$ is connected.  Clearly,
$s_\Gamma=s$. We have   $\chi_- (\Gamma)  =\vert\vert s \vert\vert$,
since
$$\vert\vert
s
\vert\vert\leq \chi_-
(\Gamma^u)=\chi_- (\Gamma) \leq \chi_- (S)=\vert\vert s \vert\vert.\eqno{\qed}$$

\skipaline {\bf 4.5\qua Proof of Theorem 3.1}\qua    Set
$\pi=\pi_1(X)$,
$H=H_1(X)$, $G=H/\Tors H$.
If $\Delta^\sigma(\pi)=0$ then $ \vert
\vert
s\vert
\vert^\sigma=0$ and $\vert \vert s\vert
\vert^\sigma-  \delta_\sigma^1 \vert s\vert_0 \leq
0 \leq \vert\vert s \vert\vert $.  Assume from now on that
$\Delta^\sigma(\pi)\neq 0$.

Fix a splitting
$H=\Tors H \times G$ and consider the ring
homomorphism $\tilde \sigma:\Bbb  Z[H] \to \C [G]$ sending $fg$
with $f\in
\Tors H, g\in G$ to $\sigma(f) g$.
By assumption,  $\gcd  \, \tilde \sigma
(E(\pi))=\Delta^\sigma(\pi)\neq 0$ so
that $\tilde \sigma
(E(\pi))\neq
0$. Pick a representative $ \sum_{g\in G} c_g g$ of
$\Delta^\sigma(\pi)$. Pick a nonzero     $\mu=\sum_{g\in G} \mu_g g\in
\tilde \sigma
(E(\pi))$ where $\mu_g\in \C$.  We
call $s\in H^1(X )$ {\sl regular} if $ s(g)
\neq
s(g')$
for any distinct $g,g'\in G$ such that $c_gc_{g'}\neq 0 $ or $\mu_g
\mu_{g'}\neq
0 $. (For $\rk G\geq 2$, this notion depends on the choice of $\mu$.
In the case $\rk G =1$ all nonzero $s$
are
regular).   The set of regular $s$ is the complement in
$H^1(X )$
of a
finite
set of sublattices of corank 1.
We call $s\in H^1(X )$ {\sl primitive} if its evaluation on a certain
element of $H$ equals $1$.

  Since   the   norms $\vert\vert ...
\vert\vert,
\vert \vert ...\vert
\vert^\sigma$, and $\vert ... \vert_0$  on $H^1(X; \Bbb
R)$ are
continuous and homogeneous, it suffices to prove that
$\vert\vert s \vert\vert \geq \vert \vert s\vert
\vert^\sigma-  \delta_\sigma^1 \vert s\vert_0$ for   primitive
regular $s\in H^1(X )$.
Fix a primitive regular $s\in H^1(X)$. (In the case $ \rk G=1$,
$s$  is
  any
generator of $H^1(X)=\Bbb  Z$).
Let
$\varphi :\Bbb  Z[H] \to
\C
[t^{\pm 1}]=\Lambda$ be   the composition of   $\tilde \sigma:\Bbb
Z[H]
\to \C [G]$   and the $\C$-linear ring
homomorphism
$\tilde s:\C [G] \to
\Lambda=\C [t^{\pm 1}]$ sending any $g\in G$ to $t^{s(g)}$.
    Recall the $\Lambda$-module $H_*^\varphi ( X)$   (see Sect. 4.2).

\skipaline {\bf Claim 1}\qua  {\sl$\Delta_0 ( H_1^\varphi ( X))\in \Lambda$ is
non-zero and
divisible by
$
(t-1)^{\delta}
\sum_{g\in
G}
c_g t^{s(g)}$ where   $\delta=1$ if $\sigma =1$ and $ \rk G\geq 2$ and
$\delta=0$
otherwise.}

\skipaline By the regularity of $s$, the   polynomial $
(t-1)^{\delta} \sum_{g } c_g t^{s(g)}$ is nonzero and its span equals
$ \delta+\vert
\vert
s\vert \vert^\sigma$.
Claim 1   implies that
$$  \dim_{\C} H_1^\varphi (X)=\span \Delta_0 ( H_1^\varphi ( X)) \geq
\delta +
\vert \vert s\vert \vert^\sigma =\delta_\sigma^1(1- \vert
s\vert_0)+
\vert \vert s\vert \vert^\sigma.$$ The
   inequality
$\vert\vert s \vert\vert \geq \vert \vert s\vert
\vert^\sigma-  \delta_\sigma^1 \vert s\vert_0$   follows  now from the
next
claim.

\skipaline {\bf Claim 2}\qua    $\vert\vert s \vert\vert \geq   \dim_{\C}\,
H_1^\varphi ( X)- \delta_\sigma^1$.

  \skipaline Now we prove Claims 1 and 2.

\proof[Proof of Claim 1]     Contracting recursively the
1-cells
of $X$
with
distinct endpoints we obtain a finite 2-dimensional CW-complex homotopy
equivalent to
$X$ and having only one $0$-cell. Since    $H_1^\varphi ( X)$
and   $\Delta^\sigma(\pi)$ are homotopy invariants, we can   assume  in
the
proof
of
Claim 1   that $X$ has only one 0-cell,   $x$.
   Consider
   the presentation of $\pi=\pi_1(X,x)$ determined by the cellular
structure of $X$.  The corresponding  Alexander matrix  is nothing but
the
matrix of
the
$\Bbb  Z[H]$-linear
boundary homomorphism $C_2(\hat X) \to C_1(\hat X)$ where $\hat X$ is
the
maximal
abelian covering of $X$.  Applying $ \varphi:  \Bbb  Z[H]\to \Lambda$
to the
entries of this matrix we obtain a presentation matrix of the
$\Lambda$-module
$H_1^{ \varphi}( X,x) $.  Therefore $ E_1(H_1^{ \varphi}( X,x))= \C
\cdot
\varphi (E
(\pi))\subset \Lambda$.

Recall that $\varphi =\tilde s \circ \tilde
\sigma $.
By definition of $\Delta^\sigma(\pi)$, we have  $  \tilde
\sigma (E
(\pi))\subset  \Delta^\sigma(\pi)\, \C [G]$. If $\sigma=1$ and $
\rk G\geq
2$, then
(2.a) implies  a   stronger inclusion  $  \tilde
\sigma (E
(\pi))\subset  \Delta^\sigma(\pi) J$
where   $J$ is the augmentation
ideal of $\C [G]$. Applying $\tilde s$, we obtain that $
\varphi (E
(\pi)) =(\tilde s \circ \tilde
\sigma) (E
(\pi))$
  is contained in the principal ideal generated by $(t-1)^{\delta}\tilde
s
(\Delta^\sigma(\pi))=
(t-1)^{ \delta} \sum_{g } c_g t^{s(g)}$. The regularity of $s$ implies
that
$
\varphi (E
(\pi))\neq 0$.
Hence
$\Delta_1 ( H_1^\varphi ( X,x))= \gcd
\varphi (E
(\pi))$ is non-zero and divisible by $
(t-1)^{\delta}
\sum_{g\in
G}
c_g t^{s(g)}$.  It remains to observe that
$\Delta_0 ( H_1^\varphi ( X))=\Delta_1 ( H_1^\varphi ( X,x))$.
Indeed, consider   the exact
sequence
$$0\to H_1^{ \varphi}( X) \to H_1^{ \varphi}( X,x) \to H_0^{
\varphi}(x) \to
H_0^{
\varphi}( X).  $$ Clearly, $ H_0^{ \varphi}(x)=\Lambda$ and
$H_0^{ \varphi}(X)=\Lambda/\varphi (I )\Lambda$ where $I$
is the augmentation ideal of $\Bbb
  Z [H]$.  The kernel $\varphi (I )\Lambda\subset \Lambda$ of the
inclusion
homomorphism
$H_0^{ \varphi}(x) \to H_0^{
\varphi} ( X) $ is a free $\Lambda$-module of rank 1.  Hence $H_1^{
\varphi}(
X,x)=H_1^{ \varphi}( X)\oplus \Lambda$ and  $ \Delta_0 (H_1^{ \varphi}(
X))
=\Delta_1 (H_1^{ \varphi}( X,x))$.\endproof

{\bf Proof of Claim 2}\qua  Consider the ring homomorphism
$\overline
\sigma=\aug \circ \tilde \sigma :
\Bbb
Z[H ] \to \C$ mapping $G\subset H$ to $1$ and mapping any $f\in
\Tors H$
to
$\sigma (f)\in \C^*$.  We call a cellular  set $S\subset X$  {\sl
bad
}
if $\overline \sigma$ is trivial on $H_1(S)$, i.e., if the composition
of the
inclusion homomorphism $ H_1(S) \to H $ with $\overline \sigma$ maps
$H_1(S)$
to
$1$.

By Lemma 4.4, there is a
  weighted graph $\Gamma=\cup_i (\Gamma_i, w_i) \subset X$ such that
$X\backslash
\Gamma$ is
connected,
$s_\Gamma=s$, and $\chi_-
(\Gamma)=\vert\vert s \vert\vert$. We first compute the
$\varphi$-twisted
homology $H_*^{ \varphi}({\Gamma_i})$ of a component  ${\Gamma_i}$ of
$\Gamma$.
   Observe that $s$
annihilates $H_1({\Gamma_i})$ and therefore
$\varphi\vert_{H_1({\Gamma_i})}$ is
the
composition
of
$\overline \sigma\vert_{H_1({\Gamma_i})}$ with the inclusion $\C\subset
\Lambda$.
Hence  $H_*^{ \varphi}({\Gamma_i})= \Lambda\otimes_{\C}
H_*^{\overline \sigma}({\Gamma_i})$.  If ${\Gamma_i}$ is bad then
$H_*^{\overline
\sigma}({\Gamma_i})$
is the usual untwisted homology of ${\Gamma_i}$ with complex
coefficients.
   If ${\Gamma_i}$ is not bad then  $H_0^{\overline
\sigma}({\Gamma_i})=0$ and
$\dim_{\C}\, H_1^{\overline \sigma}({\Gamma_i})= -\chi ({\Gamma_i}) =
\chi_-
({\Gamma_i})\geq 0$.

    Let $U=\Gamma \times [-1,1]$ be a closed regular
neighborhood of
$
\Gamma $ in
$X$ such that $\Gamma=\Gamma\times 0$. We can assume that the given
coorientation of
$\Gamma$ is determined by $\Gamma\times (0,1]\subset U$. Set
$N=\overline {X\backslash U}$. By our assumptions, $N $ is connected.
Clearly,
$N\cap
U=\partial N=\partial U$
contains
two copies  $\Gamma^{\pm}_i=\Gamma_i\times (\pm 1)$
of
each $\Gamma_i$.

The  Mayer-Vietoris  homology sequence of  the triple $(X=N\cup U, N,
U)$, gives
an exact sequence
$$H_1^{ \varphi}(X)\to H_0^{ \varphi}(N\cap U)\to
H_0^{ \varphi}(N ) \oplus H_0^{ \varphi}( U).$$
It
is clear that $s$ annihilates $H_1(N )$ and therefore
  $H_0^{\varphi}(N )=\Lambda\otimes_{\C}
H_0^{\overline \sigma}(N )= \Lambda^\beta$ where $\beta=1$ if $N$
is bad and
$\beta=0$
otherwise. The computations above show that $H_0^{
\varphi}(U)=H_0^{
\varphi}(\Gamma)=\Lambda^{\alpha}$ where $\alpha$ is the number of bad
components
of
$\Gamma$. Similarly, $H_0^{
\varphi}(N\cap U)=\Lambda^{2\alpha}$.
Therefore,   $$0=\rk_\Lambda \, H_1^{ \varphi}(X)\geq 2\alpha -
(\alpha+\beta)= \alpha-\beta.$$ Hence $\alpha \leq \beta=0,1$.

\skipaline {\bf Case $ \beta=1$}\qua  In this case     $N$
is bad and therefore all its boundary components $\Gamma^{\pm}_i$ are
bad. Thus,
all
the
components
of $\Gamma$ are bad. The inequality $\alpha\leq \beta=1$ implies that
$\Gamma$
is
connected.
Since the dual class $s$ is primitive, the weight of (the only
component of)
$\Gamma$ is equal to $1$. Thus $\vert\vert s\vert\vert=\chi_-
(\Gamma)$.

   Since $\Gamma$ and
its
complement in $X$ are bad, the group homomorphism $\overline
\sigma\vert_H:H\to \C^*$ is a composition of $ s:  H\to \Bbb  Z$ with a
certain
group
homomorphism $\Bbb
Z\to \C^*$.  Such a composition is trivial on $\Tors H$.  Hence
$\sigma=1$. Then  $H_1^\varphi ( X)= H_1(\tilde X;\C)$ where $\tilde
X\to
X$
is
the
infinite
cyclic covering determined by $s$.  To prove Claim 2, it suffices to
prove  the
inequality $\chi_- (\Gamma)   \geq \dim_{\C}\, H_1(\tilde X;
\C) -1
$.
Observe that
the graph $\Gamma$ lifts to a homeomorphic graph $\tilde \Gamma
\subset \tilde
X$  splitting $\tilde X$ into two connected pieces, $\tilde X_-$ and
$\tilde
X_+$.  Let $t$ be the generating deck transformation of the covering
$\tilde
X\to
X$ such that $t \tilde X_+ \subset \tilde X_+$.  The characteristic
polynomial
of
the action of $t$ on $H_1(\tilde X, \C)$ is $\Delta_0(H_1(\tilde
X; \C))=\Delta_0 (H_1^{ \varphi}( X)) \neq 0$.  Applying to any compact
subset of
$\tilde X$ a sufficiently big positive (resp.\ negative) power of $t$
we can
translate this subset into $\tilde X_+$ (resp.\ $\tilde X_-$).  This
implies
that
the inclusion homomorphisms $ H_1(\tilde X_-; \C) \to H_1(\tilde
X, \C)$
and $ H_1(\tilde X_+; \C) \to H_1(\tilde X, \C)$ are
surjective.  The
Mayer-Vietoris homology sequence for $\tilde X=\tilde X_+\cup \tilde
X_-$
implies
the surjectivity of the inclusion homomorphism $ H_1(\tilde \Gamma ;
\C)
\to
H_1(\tilde X, \C)$.  Computing the dimensions, we obtain $\chi_-
(\Gamma)+1=\dim_{\C}\,
H_1(\Gamma; \C)
\geq
\dim_{\C}\,
H_1(\tilde X; \C)$.

\skipaline {\bf Case $ \beta=0$}\qua  In this case $N$ is not bad and
$\alpha=\beta=0$ so that
$\Gamma
$
has no bad
components.   In particular, $\sigma\neq 1$.

Now we compute $H_1^\varphi (X)$. Let $\Gamma_1,..., \Gamma_n$ be the
components of
$\Gamma$ with
weights
$w_1,...,w_n$.  As we know
$H_0^{
\varphi}(N\cap U)=\Lambda^{2\alpha} =0 $.  The Mayer-Vietoris homology
sequence of
the triple $(X=N\cup U, N, U )$ yields that the inclusion homomorphism
$H_1^\varphi (N)
\to
H_1^\varphi (X)$ is surjective and its kernel is generated by the
vectors $in
(t^{w_i} x-f_i^\varphi (x))$ where $i=1,..., n$; $x$ runs over $
H_1^\varphi(\Gamma^{+}_i)$;
$f_i^\varphi:
H_1^\varphi(\Gamma^{+}_i)\to H_1^\varphi(\Gamma^{-}_i)$ is the
isomorphism
induced
by the
natural homeomorphisms $\Gamma^{+}_i\approx \Gamma_i\approx
\Gamma^{-}_i$; and
$in$
is the
inclusion homomorphism $H_1^\varphi (\partial N) \to H_1^\varphi (N)$.
We
claim
that $in $ is surjective.  Indeed, since $s$ annihilates $H_1(N)$ we
have
$H_1^{
\varphi}(N)= \Lambda\otimes_{\C} H_1^{\overline \sigma}(N)$ and
$H_1^\varphi
(\partial N) = \Lambda\otimes_{\C} H_1^{\overline \sigma}(\partial
N)$.
Moreover, $in =\id_\Lambda \otimes_{\C}\, j$ where $j:
H_1^{\overline
\sigma}(\partial N) \to H_1^{\overline \sigma} ( N)$ is the inclusion
homomorphism.
If $j$ is not surjective then the cokernel of $ in$ is a free
$\Lambda$-module
of
rank $\geq 1$.  On the other hand, this cokernel
is
a
quotient of the finite dimensional $\C$-linear space $
H_1^{\varphi}(X)$.
This
contradiction shows that both $j$ and $in $ must be surjective.
Therefore
$H_1^\varphi (X)$ is the quotient of $$H_1^{\varphi}(\partial
N)=\bigoplus_{i=1}^n
(H_1^{\varphi}(\Gamma^{+}_i)\oplus H_1^{\varphi}(\Gamma^{-}_i))
=\Lambda
\otimes_{\C}
\bigoplus_{i=1}^n (H_1^{\overline \sigma}(\Gamma^{+}_i)\oplus
H_1^{\overline
\sigma}(\Gamma^{-}_i)) $$ by $\Ker\, in=\Lambda \otimes_{\C} \Ker
\,j $ and
the
vectors
$
t^{w_i} x- f_i^{\overline\sigma}  (x) $ where $i$ runs over $1,..., n$;
$x$ runs
over $
H_1^{\overline\sigma} (\Gamma^{+}_i)$;
$f_i^{\overline\sigma}:
H_1^{\overline\sigma} (\Gamma^{+}_i)\to H_1^{\overline\sigma}
(\Gamma^{-}_i)$
is the
isomorphism
induced
by the
   homeomorphism  $\Gamma^{+}_i\approx   \Gamma^{-}_i$.  Hence
$H_1^\varphi (X)$ is a
quotient of $
\Lambda \otimes_{\C} (\oplus_{i} H_1^{\overline
\sigma}(\Gamma^{+}_i) )$
by
   vectors of type $x_1+...+x_n +t^{w_1} y_1+...+ t^{w_n}
y_n
$
where $x_i,y_i\in H_1^{\overline \sigma}(\Gamma^{+}_i)$ for all $i $.
Consider
the
corresponding presentation matrix of $H_1^\varphi (X)$ with respect to
certain
bases in the $\C$-linear spaces $H_1^{\overline
\sigma}(\Gamma^{+}_1),...,
H_1^{\overline \sigma}(\Gamma^{+}_n) $.  The column of this matrix
corresponding
to
any
basis vector in $H_1^{\overline \sigma}(\Gamma^{+}_i) $ has entries of
type $a+
t^{w_i}b$ with $a,b\in \C$.  Therefore the ideal $E_0(H_1^\varphi
(X))$ is
generated by
Laurent polynomials whose span does not exceed
$$\sum_{i=1}^n w_i\, \dim_{\C}
\, H_1^{\overline \sigma}(\Gamma^{+}_i)   = \sum_{i=1}^n w_i \, \chi_-
(\Gamma^+_i)=\chi_-
(\Gamma)= \vert \vert s\vert\vert .$$
Therefore
$$\dim_{\C}\,
H_1^\varphi ( X)- \delta_\sigma^1=\dim_{\C}\, H_1^\varphi ( X)=
{\rm {span}}\, \Delta_0 (H_1^\varphi (X))\leq   \vert \vert
s\vert\vert .\eqno{\qed}$$

\Addressesr
\end{document}

%% file: agtout.tex
\def\ifplaintex{\expandafter\ifx\csname documentclass\endcsname\relax}

\def\gtp{{\mathsurround=0pt\it $\cal G\mskip-2mu$eometry \&\ 
$\cal T\!\!$opology $\cal P\!$ublications}}  

\def\Addressesr{\bigskip
{\small \parskip 0pt \leftskip 0pt \rightskip 0pt plus 1fil \def\\{\par}
\sl\theaddress\par
\medskip
\rm Email:\stdspace\tt\theemail\hfill\rm Received:\qua\receiveddate \par}}

\def\recd{{\small Received:\qua\receiveddate\ifx\reviseddate\relax
\else\qquad Revised:\qua\reviseddate\fi\par}} 

\def\lognumber#1{\def\thelognumber{#1}}
\def\volumenumber#1{\def\thevolumenumber{#1}}
\def\volumeyear#1{\def\thevolumeyear{#1}}
\def\papernumber#1{\def\thepapernumber{#1}}
\def\pagenumbers#1#2{\def\startpage{#1}\def\finishpage{#2}}
\def\published#1{\def\publishdate{#1}}

\def\received#1{\def\receiveddate{#1}}

\def\accepted#1{\def\accepteddate{#1}}

\def\asciiaddress#1{\def\theasciiaddress{#1}}

\long\def\asciiabstract#1{\long\def\theasciiabstract{#1}}

\let\\\par\let\thelognumber\relax\let\thevolumenumber\relax
\let\thepapernumber\relax\let\thevolumeyear\relax\let\startpage\relax
\let\finishpage\relax\let\publishdate\relax\let\receiveddate\relax
\let\reviseddate\relax\let\accepteddate\relax\let\theasciititle\relax
\let\theasciiauthors\relax\let\theasciiaddress\relax
\let\theasciiabstract\relax

\let\theasciiemail\relax

\ifplaintex
\font\logobig=cmssbx10 scaled 3836
\font\logomed=cmssbx10 scaled 2557
\else
\font\logobig=cmssbx10 scaled 4200
\font\logomed=cmssbx10 scaled 2800
\fi

\long\def\makeagttitle{   
\count0=\startpage
\agt\hfill      
\hbox to 45truept{\vbox to 0pt{\vglue -13truept{\logomed A\kern -.37em{\logobig 
T}\kern -.38em G}\vss}\hss}
\break
{\small Volume \thevolumenumber\ (\thevolumeyear)
\startpage--\finishpage\nl
Published: \publishdate}

\vglue .25truein

{\parskip=0pt\leftskip 0pt plus
1fil\def\\{\par\smallskip}{\Large\bf\thetitle}\par\medskip} \vglue
0.05truein

{\parskip=0pt\leftskip 0pt plus 1fil\def\\{\par}{\sc\theauthors}
\par\medskip}%
 
\vglue 0.03truein 

{\small\leftskip 25truept\rightskip 25truept{\bf Abstract}\stdspace\theabstract

{\bf AMS Classification}\stdspace\theprimaryclass
\ifx\thesecondaryclass\relax\else; \thesecondaryclass\fi\par
{\bf Keywords}\stdspace \thekeywords\par}\vglue 7truept

}   

\ifplaintex
\font\phead=cmsl9 scaled 950
\font\pnum=cmbx10 scaled 913
\font\pfoot=cmsl9 scaled 950
\headline{\vbox to 0pt{\vskip -4.5mm\line{\small\phead\ifnum
\count0=\startpage ISSN 1472-2739 (on-line) 1472-2747 (printed)
\hfill {\pnum\folio}\else\ifodd\count0\def\\{ }%
\ifx\theshorttitle\relax\thetitle\else\theshorttitle\fi\hfill{\pnum\folio}
\else\def\\{ and }{\pnum\folio}\hfill\ifx\theshortauthors\relax\theauthors
\else\theshortauthors\fi\fi\fi}\vss}}
\footline{\vbox to 0pt{\vglue 0mm\line{\small\pfoot\ifnum\count0=\startpage
\copyright\ \gtp\hfill\else
\agt, Volume \thevolumenumber\ (\thevolumeyear)\hfill\fi}\vss}}
\else
\font=cmsl9 scaled 1050
\font\lnum=cmbx10 
\font=cmsl9 scaled 1050
\makeatletter
\def\@oddhead{{\small\lhead\ifnum\count0=\startpage ISSN 1472-2739 
(on-line) 1472-2747 (printed)\hfill {\lnum\number\count0}\else\ifodd\count0
\def\\{ }\ifx\theshorttitle\relax \thetitle \else\theshorttitle\fi\hfill
{\lnum\number\count0}\else\def\\{ and }{\lnum\number\count0}
\hfill\ifx\theshortauthors\relax 
\theauthors\else\theshortauthors\fi\fi\fi}}\def\@evenhead{\@oddhead}
\def\@oddfoot{\small\lfoot\ifnum\count0=\startpage\copyright\ \gtp\hfill\else
\agt, Volume \thevolumenumber\ (\thevolumeyear)\hfill\fi}
\def\@evenfoot{\@oddfoot}
\makeatother
\fi
\let\maketitlepage\makeagttitle
\let\makeshorttitle\maketitlepage
\let\maketitle\maketitlepage

\newwrite\gtoutfile
\long\gdef\makeheadfile{  
{\def\\{, }\def\s{ }
\immediate\openout\gtoutfile head.xxx
\immediate\write\gtoutfile{To: math@arxiv.org}
\immediate\write\gtoutfile{Subject: put OR rep NNNNN:ppppp}
\immediate\write\gtoutfile{--text follows this line--}
\immediate\write\gtoutfile{Proxy-for: \ifx\theasciiauthors\relax
\theauthors\else\theasciiauthors\fi\s<\ifx\theasciiemail\relax\theemail\else\theasciiemail\fi>}
\immediate\write\gtoutfile{\noexpand\\}
\immediate\write\gtoutfile{Authors: \ifx\theasciiauthors\relax
\theauthors\else\theasciiauthors\fi}
{\def\\{ }\immediate\write\gtoutfile{Title: \ifx\theasciititle\relax
\thetitle\else\theasciititle\fi}}
\immediate\write\gtoutfile{Subj-class: GT or SG, GR etc}
\immediate\write\gtoutfile{MSC-class: \theprimaryclass\ifx\thesecondaryclass\relax\else, \thesecondaryclass\fi}
\immediate\write\gtoutfile{Journal-ref: Algebr. Geom. Topol. \thevolumenumber\s
(\thevolumeyear) \startpage-\finishpage}
\immediate\write\gtoutfile{Comments: Published by Algebraic and
Geometric Topology at}
\immediate\write\gtoutfile{\s\s\s  http://www.maths.warwick.ac.uk/agt/AGTVol\thevolumenumber/agt-\thevolumenumber-\thepapernumber.abs.html}
\immediate\write\gtoutfile{\noexpand\\}
\immediate\write\gtoutfile{}
\ifx\theasciiabstract\relax
\immediate\write\gtoutfile{\theabstract}\else
\immediate\write\gtoutfile{\theasciiabstract}\fi
\immediate\write\gtoutfile{}
\immediate\write\gtoutfile{\noexpand\\}
\immediate\write\gtoutfile{}
\immediate\closeout\gtoutfile}}  

\def\maketitlepage{\makeagttitle\makeheadfile}
\let\makeshorttitle\maketitlepage
\let\maketitle\maketitlepage

\def\ifplaintex{\expandafter\ifx\csname documentclass\endcsname\relax}

\def\gtp{{\mathsurround=0pt\it $\cal G\mskip-2mu$eometry \&\ 
$\cal T\!\!$opology $\cal P\!$ublications}}  

\def\Addressesr{\bigskip
{\small \parskip 0pt \leftskip 0pt \rightskip 0pt plus 1fil \def\\{\par}
\sl\theaddress\par
\medskip
\rm Email:\stdspace\tt\theemail\hfill\rm Received:\qua\receiveddate \par}}

\def\recd{{\small Received:\qua\receiveddate\ifx\reviseddate\relax
\else\qquad Revised:\qua\reviseddate\fi\par}} 

\def\lognumber#1{\def\thelognumber{#1}}
\def\volumenumber#1{\def\thevolumenumber{#1}}
\def\volumeyear#1{\def\thevolumeyear{#1}}
\def\papernumber#1{\def\thepapernumber{#1}}
\def\pagenumbers#1#2{\def\startpage{#1}\def\finishpage{#2}}
\def\published#1{\def\publishdate{#1}}

\def\received#1{\def\receiveddate{#1}}

\def\accepted#1{\def\accepteddate{#1}}

\def\asciiaddress#1{\def\theasciiaddress{#1}}

\long\def\asciiabstract#1{\long\def\theasciiabstract{#1}}

\let\\\par\let\thelognumber\relax\let\thevolumenumber\relax
\let\thepapernumber\relax\let\thevolumeyear\relax\let\startpage\relax
\let\finishpage\relax\let\publishdate\relax\let\receiveddate\relax
\let\reviseddate\relax\let\accepteddate\relax\let\theasciititle\relax
\let\theasciiauthors\relax\let\theasciiaddress\relax
\let\theasciiabstract\relax

\let\theasciiemail\relax

\ifplaintex
\font\logobig=cmssbx10 scaled 3836
\font\logomed=cmssbx10 scaled 2557
\else
\font\logobig=cmssbx10 scaled 4200
\font\logomed=cmssbx10 scaled 2800
\fi

\long\def\makeagttitle{   
\count0=\startpage
\agt\hfill      
\hbox to 45truept{\vbox to 0pt{\vglue -13truept{\logomed A\kern -.37em{\logobig 
T}\kern -.38em G}\vss}\hss}
\break
{\small Volume \thevolumenumber\ (\thevolumeyear)
\startpage--\finishpage\nl
Published: \publishdate}

\vglue .25truein

{\parskip=0pt\leftskip 0pt plus
1fil\def\\{\par\smallskip}{\Large\bf\thetitle}\par\medskip} \vglue
0.05truein

{\parskip=0pt\leftskip 0pt plus 1fil\def\\{\par}{\sc\theauthors}
\par\medskip}%
 
\vglue 0.03truein 

{\small\leftskip 25truept\rightskip 25truept{\bf Abstract}\stdspace\theabstract

{\bf AMS Classification}\stdspace\theprimaryclass
\ifx\thesecondaryclass\relax\else; \thesecondaryclass\fi\par
{\bf Keywords}\stdspace \thekeywords\par}\vglue 7truept

}   

\ifplaintex
\font\phead=cmsl9 scaled 950
\font\pnum=cmbx10 scaled 913
\font\pfoot=cmsl9 scaled 950
\headline{\vbox to 0pt{\vskip -4.5mm\line{\small\phead\ifnum
\count0=\startpage ISSN 1472-2739 (on-line) 1472-2747 (printed)
\hfill {\pnum\folio}\else\ifodd\count0\def\\{ }%
\ifx\theshorttitle\relax\thetitle\else\theshorttitle\fi\hfill{\pnum\folio}
\else\def\\{ and }{\pnum\folio}\hfill\ifx\theshortauthors\relax\theauthors
\else\theshortauthors\fi\fi\fi}\vss}}
\footline{\vbox to 0pt{\vglue 0mm\line{\small\pfoot\ifnum\count0=\startpage
\copyright\ \gtp\hfill\else
\agt, Volume \thevolumenumber\ (\thevolumeyear)\hfill\fi}\vss}}
\else
\font=cmsl9 scaled 1050
\font\lnum=cmbx10 
\font=cmsl9 scaled 1050
\makeatletter
\def\@oddhead{{\small\lhead\ifnum\count0=\startpage ISSN 1472-2739 
(on-line) 1472-2747 (printed)\hfill {\lnum\number\count0}\else\ifodd\count0
\def\\{ }\ifx\theshorttitle\relax \thetitle \else\theshorttitle\fi\hfill
{\lnum\number\count0}\else\def\\{ and }{\lnum\number\count0}
\hfill\ifx\theshortauthors\relax 
\theauthors\else\theshortauthors\fi\fi\fi}}\def\@evenhead{\@oddhead}
\def\@oddfoot{\small\lfoot\ifnum\count0=\startpage\copyright\ \gtp\hfill\else
\agt, Volume \thevolumenumber\ (\thevolumeyear)\hfill\fi}
\def\@evenfoot{\@oddfoot}
\makeatother
\fi
\let\maketitlepage\makeagttitle
\let\makeshorttitle\maketitlepage
\let\maketitle\maketitlepage

\newwrite\gtoutfile
\long\gdef\makeheadfile{  
{\def\\{, }\def\s{ }
\immediate\openout\gtoutfile head.xxx
\immediate\write\gtoutfile{To: math@arxiv.org}
\immediate\write\gtoutfile{Subject: put OR rep NNNNN:ppppp}
\immediate\write\gtoutfile{--text follows this line--}
\immediate\write\gtoutfile{Proxy-for: \ifx\theasciiauthors\relax
\theauthors\else\theasciiauthors\fi\s<\ifx\theasciiemail\relax\theemail\else\theasciiemail\fi>}
\immediate\write\gtoutfile{\noexpand\\}
\immediate\write\gtoutfile{Authors: \ifx\theasciiauthors\relax
\theauthors\else\theasciiauthors\fi}
{\def\\{ }\immediate\write\gtoutfile{Title: \ifx\theasciititle\relax
\thetitle\else\theasciititle\fi}}
\immediate\write\gtoutfile{Subj-class: GT or SG, GR etc}
\immediate\write\gtoutfile{MSC-class: \theprimaryclass\ifx\thesecondaryclass\relax\else, \thesecondaryclass\fi}
\immediate\write\gtoutfile{Journal-ref: Algebr. Geom. Topol. \thevolumenumber\s
(\thevolumeyear) \startpage-\finishpage}
\immediate\write\gtoutfile{Comments: Published by Algebraic and
Geometric Topology at}
\immediate\write\gtoutfile{\s\s\s  http://www.maths.warwick.ac.uk/agt/AGTVol\thevolumenumber/agt-\thevolumenumber-\thepapernumber.abs.html}
\immediate\write\gtoutfile{\noexpand\\}
\immediate\write\gtoutfile{}
\ifx\theasciiabstract\relax
\immediate\write\gtoutfile{\theabstract}\else
\immediate\write\gtoutfile{\theasciiabstract}\fi
\immediate\write\gtoutfile{}
\immediate\write\gtoutfile{\noexpand\\}
\immediate\write\gtoutfile{}
\immediate\closeout\gtoutfile}}  

\def\maketitlepage{\makeagttitle\makeheadfile}
\let\makeshorttitle\maketitlepage
\let\maketitle\maketitlepage

\def\ifplaintex{\expandafter\ifx\csname documentclass\endcsname\relax}

\def\gtp{{\mathsurround=0pt\it $\cal G\mskip-2mu$eometry \&\ 
$\cal T\!\!$opology $\cal P\!$ublications}}  

\def\Addressesr{\bigskip
{\small \parskip 0pt \leftskip 0pt \rightskip 0pt plus 1fil \def\\{\par}
\sl\theaddress\par
\medskip
\rm Email:\stdspace\tt\theemail\hfill\rm Received:\qua\receiveddate \par}}

\def\recd{{\small Received:\qua\receiveddate\ifx\reviseddate\relax
\else\qquad Revised:\qua\reviseddate\fi\par}} 

\def\lognumber#1{\def\thelognumber{#1}}
\def\volumenumber#1{\def\thevolumenumber{#1}}
\def\volumeyear#1{\def\thevolumeyear{#1}}
\def\papernumber#1{\def\thepapernumber{#1}}
\def\pagenumbers#1#2{\def\startpage{#1}\def\finishpage{#2}}
\def\published#1{\def\publishdate{#1}}

\def\received#1{\def\receiveddate{#1}}

\def\accepted#1{\def\accepteddate{#1}}

\def\asciiaddress#1{\def\theasciiaddress{#1}}

\long\def\asciiabstract#1{\long\def\theasciiabstract{#1}}

\let\\\par\let\thelognumber\relax\let\thevolumenumber\relax
\let\thepapernumber\relax\let\thevolumeyear\relax\let\startpage\relax
\let\finishpage\relax\let\publishdate\relax\let\receiveddate\relax
\let\reviseddate\relax\let\accepteddate\relax\let\theasciititle\relax
\let\theasciiauthors\relax\let\theasciiaddress\relax
\let\theasciiabstract\relax

\let\theasciiemail\relax

\ifplaintex
\font\logobig=cmssbx10 scaled 3836
\font\logomed=cmssbx10 scaled 2557
\else
\font\logobig=cmssbx10 scaled 4200
\font\logomed=cmssbx10 scaled 2800
\fi

\long\def\makeagttitle{   
\count0=\startpage
\agt\hfill      
\hbox to 45truept{\vbox to 0pt{\vglue -13truept{\logomed A\kern -.37em{\logobig 
T}\kern -.38em G}\vss}\hss}
\break
{\small Volume \thevolumenumber\ (\thevolumeyear)
\startpage--\finishpage\nl
Published: \publishdate}

\vglue .25truein

{\parskip=0pt\leftskip 0pt plus
1fil\def\\{\par\smallskip}{\Large\bf\thetitle}\par\medskip} \vglue
0.05truein

{\parskip=0pt\leftskip 0pt plus 1fil\def\\{\par}{\sc\theauthors}
\par\medskip}%
 
\vglue 0.03truein 

{\small\leftskip 25truept\rightskip 25truept{\bf Abstract}\stdspace\theabstract

{\bf AMS Classification}\stdspace\theprimaryclass
\ifx\thesecondaryclass\relax\else; \thesecondaryclass\fi\par
{\bf Keywords}\stdspace \thekeywords\par}\vglue 7truept

}   

\ifplaintex
\font\phead=cmsl9 scaled 950
\font\pnum=cmbx10 scaled 913
\font\pfoot=cmsl9 scaled 950
\headline{\vbox to 0pt{\vskip -4.5mm\line{\small\phead\ifnum
\count0=\startpage ISSN 1472-2739 (on-line) 1472-2747 (printed)
\hfill {\pnum\folio}\else\ifodd\count0\def\\{ }%
\ifx\theshorttitle\relax\thetitle\else\theshorttitle\fi\hfill{\pnum\folio}
\else\def\\{ and }{\pnum\folio}\hfill\ifx\theshortauthors\relax\theauthors
\else\theshortauthors\fi\fi\fi}\vss}}
\footline{\vbox to 0pt{\vglue 0mm\line{\small\pfoot\ifnum\count0=\startpage
\copyright\ \gtp\hfill\else
\agt, Volume \thevolumenumber\ (\thevolumeyear)\hfill\fi}\vss}}
\else
\font=cmsl9 scaled 1050
\font\lnum=cmbx10 
\font=cmsl9 scaled 1050
\makeatletter
\def\@oddhead{{\small\lhead\ifnum\count0=\startpage ISSN 1472-2739 
(on-line) 1472-2747 (printed)\hfill {\lnum\number\count0}\else\ifodd\count0
\def\\{ }\ifx\theshorttitle\relax \thetitle \else\theshorttitle\fi\hfill
{\lnum\number\count0}\else\def\\{ and }{\lnum\number\count0}
\hfill\ifx\theshortauthors\relax 
\theauthors\else\theshortauthors\fi\fi\fi}}\def\@evenhead{\@oddhead}
\def\@oddfoot{\small\lfoot\ifnum\count0=\startpage\copyright\ \gtp\hfill\else
\agt, Volume \thevolumenumber\ (\thevolumeyear)\hfill\fi}
\def\@evenfoot{\@oddfoot}
\makeatother
\fi
\let\maketitlepage\makeagttitle
\let\makeshorttitle\maketitlepage
\let\maketitle\maketitlepage

\newwrite\gtoutfile
\long\gdef\makeheadfile{  
{\def\\{, }\def\s{ }
\immediate\openout\gtoutfile head.xxx
\immediate\write\gtoutfile{To: math@arxiv.org}
\immediate\write\gtoutfile{Subject: put OR rep NNNNN:ppppp}
\immediate\write\gtoutfile{--text follows this line--}
\immediate\write\gtoutfile{Proxy-for: \ifx\theasciiauthors\relax
\theauthors\else\theasciiauthors\fi\s<\ifx\theasciiemail\relax\theemail\else\theasciiemail\fi>}
\immediate\write\gtoutfile{\noexpand\\}
\immediate\write\gtoutfile{Authors: \ifx\theasciiauthors\relax
\theauthors\else\theasciiauthors\fi}
{\def\\{ }\immediate\write\gtoutfile{Title: \ifx\theasciititle\relax
\thetitle\else\theasciititle\fi}}
\immediate\write\gtoutfile{Subj-class: GT or SG, GR etc}
\immediate\write\gtoutfile{MSC-class: \theprimaryclass\ifx\thesecondaryclass\relax\else, \thesecondaryclass\fi}
\immediate\write\gtoutfile{Journal-ref: Algebr. Geom. Topol. \thevolumenumber\s
(\thevolumeyear) \startpage-\finishpage}
\immediate\write\gtoutfile{Comments: Published by Algebraic and
Geometric Topology at}
\immediate\write\gtoutfile{\s\s\s  http://www.maths.warwick.ac.uk/agt/AGTVol\thevolumenumber/agt-\thevolumenumber-\thepapernumber.abs.html}
\immediate\write\gtoutfile{\noexpand\\}
\immediate\write\gtoutfile{}
\ifx\theasciiabstract\relax
\immediate\write\gtoutfile{\theabstract}\else
\immediate\write\gtoutfile{\theasciiabstract}\fi
\immediate\write\gtoutfile{}
\immediate\write\gtoutfile{\noexpand\\}
\immediate\write\gtoutfile{}
\immediate\closeout\gtoutfile}}  

\def\maketitlepage{\makeagttitle\makeheadfile}
\let\makeshorttitle\maketitlepage
\let\maketitle\maketitlepage

\def\ifplaintex{\expandafter\ifx\csname documentclass\endcsname\relax}

\def\gtp{{\mathsurround=0pt\it $\cal G\mskip-2mu$eometry \&\ 
$\cal T\!\!$opology $\cal P\!$ublications}}  

\def\Addressesr{\bigskip
{\small \parskip 0pt \leftskip 0pt \rightskip 0pt plus 1fil \def\\{\par}
\sl\theaddress\par
\medskip
\rm Email:\stdspace\tt\theemail\hfill\rm Received:\qua\receiveddate \par}}

\def\recd{{\small Received:\qua\receiveddate\ifx\reviseddate\relax
\else\qquad Revised:\qua\reviseddate\fi\par}} 

\def\lognumber#1{\def\thelognumber{#1}}
\def\volumenumber#1{\def\thevolumenumber{#1}}
\def\volumeyear#1{\def\thevolumeyear{#1}}
\def\papernumber#1{\def\thepapernumber{#1}}
\def\pagenumbers#1#2{\def\startpage{#1}\def\finishpage{#2}}
\def\published#1{\def\publishdate{#1}}

\def\received#1{\def\receiveddate{#1}}

\def\accepted#1{\def\accepteddate{#1}}

\def\asciiaddress#1{\def\theasciiaddress{#1}}

\long\def\asciiabstract#1{\long\def\theasciiabstract{#1}}

\let\\\par\let\thelognumber\relax\let\thevolumenumber\relax
\let\thepapernumber\relax\let\thevolumeyear\relax\let\startpage\relax
\let\finishpage\relax\let\publishdate\relax\let\receiveddate\relax
\let\reviseddate\relax\let\accepteddate\relax\let\theasciititle\relax
\let\theasciiauthors\relax\let\theasciiaddress\relax
\let\theasciiabstract\relax

\let\theasciiemail\relax

\ifplaintex
\font\logobig=cmssbx10 scaled 3836
\font\logomed=cmssbx10 scaled 2557
\else
\font\logobig=cmssbx10 scaled 4200
\font\logomed=cmssbx10 scaled 2800
\fi

\long\def\makeagttitle{   
\count0=\startpage
\agt\hfill      
\hbox to 45truept{\vbox to 0pt{\vglue -13truept{\logomed A\kern -.37em{\logobig 
T}\kern -.38em G}\vss}\hss}
\break
{\small Volume \thevolumenumber\ (\thevolumeyear)
\startpage--\finishpage\nl
Published: \publishdate}

\vglue .25truein

{\parskip=0pt\leftskip 0pt plus
1fil\def\\{\par\smallskip}{\Large\bf\thetitle}\par\medskip} \vglue
0.05truein

{\parskip=0pt\leftskip 0pt plus 1fil\def\\{\par}{\sc\theauthors}
\par\medskip}%
 
\vglue 0.03truein 

{\small\leftskip 25truept\rightskip 25truept{\bf Abstract}\stdspace\theabstract

{\bf AMS Classification}\stdspace\theprimaryclass
\ifx\thesecondaryclass\relax\else; \thesecondaryclass\fi\par
{\bf Keywords}\stdspace \thekeywords\par}\vglue 7truept

}   

\ifplaintex
\font\phead=cmsl9 scaled 950
\font\pnum=cmbx10 scaled 913
\font\pfoot=cmsl9 scaled 950
\headline{\vbox to 0pt{\vskip -4.5mm\line{\small\phead\ifnum
\count0=\startpage ISSN 1472-2739 (on-line) 1472-2747 (printed)
\hfill {\pnum\folio}\else\ifodd\count0\def\\{ }%
\ifx\theshorttitle\relax\thetitle\else\theshorttitle\fi\hfill{\pnum\folio}
\else\def\\{ and }{\pnum\folio}\hfill\ifx\theshortauthors\relax\theauthors
\else\theshortauthors\fi\fi\fi}\vss}}
\footline{\vbox to 0pt{\vglue 0mm\line{\small\pfoot\ifnum\count0=\startpage
\copyright\ \gtp\hfill\else
\agt, Volume \thevolumenumber\ (\thevolumeyear)\hfill\fi}\vss}}
\else
\font=cmsl9 scaled 1050
\font\lnum=cmbx10 
\font=cmsl9 scaled 1050
\makeatletter
\def\@oddhead{{\small\lhead\ifnum\count0=\startpage ISSN 1472-2739 
(on-line) 1472-2747 (printed)\hfill {\lnum\number\count0}\else\ifodd\count0
\def\\{ }\ifx\theshorttitle\relax \thetitle \else\theshorttitle\fi\hfill
{\lnum\number\count0}\else\def\\{ and }{\lnum\number\count0}
\hfill\ifx\theshortauthors\relax 
\theauthors\else\theshortauthors\fi\fi\fi}}\def\@evenhead{\@oddhead}
\def\@oddfoot{\small\lfoot\ifnum\count0=\startpage\copyright\ \gtp\hfill\else
\agt, Volume \thevolumenumber\ (\thevolumeyear)\hfill\fi}
\def\@evenfoot{\@oddfoot}
\makeatother
\fi
\let\maketitlepage\makeagttitle
\let\makeshorttitle\maketitlepage
\let\maketitle\maketitlepage

\newwrite\gtoutfile
\long\gdef\makeheadfile{  
{\def\\{, }\def\s{ }
\immediate\openout\gtoutfile head.xxx
\immediate\write\gtoutfile{To: math@arxiv.org}
\immediate\write\gtoutfile{Subject: put OR rep NNNNN:ppppp}
\immediate\write\gtoutfile{--text follows this line--}
\immediate\write\gtoutfile{Proxy-for: \ifx\theasciiauthors\relax
\theauthors\else\theasciiauthors\fi\s<\ifx\theasciiemail\relax\theemail\else\theasciiemail\fi>}
\immediate\write\gtoutfile{\noexpand\\}
\immediate\write\gtoutfile{Authors: \ifx\theasciiauthors\relax
\theauthors\else\theasciiauthors\fi}
{\def\\{ }\immediate\write\gtoutfile{Title: \ifx\theasciititle\relax
\thetitle\else\theasciititle\fi}}
\immediate\write\gtoutfile{Subj-class: GT or SG, GR etc}
\immediate\write\gtoutfile{MSC-class: \theprimaryclass\ifx\thesecondaryclass\relax\else, \thesecondaryclass\fi}
\immediate\write\gtoutfile{Journal-ref: Algebr. Geom. Topol. \thevolumenumber\s
(\thevolumeyear) \startpage-\finishpage}
\immediate\write\gtoutfile{Comments: Published by Algebraic and
Geometric Topology at}
\immediate\write\gtoutfile{\s\s\s  http://www.maths.warwick.ac.uk/agt/AGTVol\thevolumenumber/agt-\thevolumenumber-\thepapernumber.abs.html}
\immediate\write\gtoutfile{\noexpand\\}
\immediate\write\gtoutfile{}
\ifx\theasciiabstract\relax
\immediate\write\gtoutfile{\theabstract}\else
\immediate\write\gtoutfile{\theasciiabstract}\fi
\immediate\write\gtoutfile{}
\immediate\write\gtoutfile{\noexpand\\}
\immediate\write\gtoutfile{}
\immediate\closeout\gtoutfile}}  

\def\maketitlepage{\makeagttitle\makeheadfile}
\let\makeshorttitle\maketitlepage
\let\maketitle\maketitlepage

%% file: agt-2-7.bbl
\begin{thebibliography}



  \bibitem{Au} {\bf D Auckly}, {\it The Thurston norm and
three-dimensional
Seiberg-Witten
theory}, Osaka J.  Math.  33 (1996)  737--750.
 

  \bibitem{Fox} {\bf R  Fox},  {\it Free differential calculus.  II.  The
isomorphism  problem
of
groups}, Ann.  of Math.  (2) 59 (1954)  196--210.

  \bibitem{Kr} {\bf P Kronheimer}, {\it Minimal genus in $S\sp 1\times M\sp
3$},  Invent.
Math.
135 (1999)  45--61.

\bibitem{KM1} {\bf P  Kronheimer}, {\bf T  Mrowka},
  {\it The genus of embedded surfaces
in the
projective plane}, Math.  Res.  Lett.  1 (1994)  797--808.

\bibitem{KM2} {\bf P  Kronheimer}, {\bf T  Mrowka}, {\it Scalar
curvature and the
Thurston norm},
Math.  Res.  Lett.  4 (1997) 931--937.

 

\bibitem{McM} {\bf C McMullen},    {\sl The Alexander polynomial of a
3-manifold and
the
Thur\-ston norm on cohomology}, Preprint (1998).


\bibitem{Th} {\bf  W   Thurston},
   {\it A norm for the homology of $3$-manifolds}, Mem.
Amer.
Math.  Soc.  59 (1986), no.  339,   99--130.

 
\bibitem{Tu1} {\bf  V Turaev}, {\it Torsion invariants of
$Spin^c$-structures on
3-manifolds},
Math.  Research Letters 4 (1997)  679--695.


 
  \bibitem{Tu2} {\bf  V Turaev}, {\it Introduction to Combinatorial
Torsions.  Notes
taken
by
Felix Schlenk}, Lectures in Mathematics ETH Z\"urich.  Birkh\"auser
Verlag,
Basel (2001).

\end{thebibliography}
